\documentclass[onecolumn]{article}
\usepackage{setspace}
\usepackage{graphicx}    %  
\usepackage{epstopdf}
\usepackage{amssymb}
\usepackage{amsmath}
\usepackage{dsfont}
\usepackage{bm}  %

\usepackage{lscape}

\usepackage{color}

\newcommand{\eod}{{\hfill $\blacksquare$}}

\newcommand{\dx}{\mathrm{d}x}

\newcommand{\dB}{\mathrm{d}B}

\newcommand{\dif}{\mathrm{d}}

\usepackage{picins}
\usepackage{multirow}
\newtheorem{prop}{Proposition}

\newtheorem{lemma}{Lemma}
\newtheorem{remark}{Remark}

\usepackage{balance}
\allowdisplaybreaks

\makeatletter
\def\blfootnote{\xdef\@thefnmark{}\@footnotetext}
\makeatother

\usepackage{keyval}
\makeatletter
\define@key{cascading}{width}{\cascading@wd=#1}
\define@key{cascading}{sep}{\def\cascading@sep{#1}}
\newdimen\cascading@wd
\newcommand{\cascadingblocks}[2][]{%
  \setkeys{cascading}{sep=2ex,#1}%
  \leavevmode\vbox{\offinterlineskip
    \@for\next:=#2\do{%
      \settowidth{\dimen@}{\next}%
      \ifdim\dimen@>\cascading@wd
        \cascading@wd=\dimen@
      \fi
    }%
    \@for\next:=#2\do{%
      \cascading@rule
      \hbox{\fbox{\hb@xt@\cascading@wd{\hss\next\hss}}}%
    }
  }
}
\def\cascading@rule{%
  \def\cascading@rule{%
    \hb@xt@\dimexpr\cascading@wd+2\fboxsep+2\fboxrule\relax
      {\hss\vrule\@height\cascading@sep\hss}%
  }%
}
\makeatother

\usepackage{smartdiagram}
\usepackage{multirow}
\def\ind{{\mathchoice {\rm 1\mskip-4mu l} {\rm 1\mskip-4mu l}
{\rm 1\mskip-4.5mu l} {\rm 1\mskip-5mu l}}}
\begin{document}
%\begin{landscape}

\title{Semi-Explicit Solutions to some Non-Linear Non-Quadratic Mean-Field-Type Games: \\  A Direct Method}
\author{Julian Barreiro-Gomez,  Tyrone E. Duncan, \\ Bozenna Pasik-Duncan, and Hamidou Tembine 
\thanks{Julian Barreiro-Gomez and Hamidou Tembine are with Learning \& Game Theory Laboratory,  New York University Abu Dhabi, (e-mails: jbarreiro@nyu.edu, tembine@nyu.edu).}
\thanks{Tyrone Duncan and Bozenna Pasik-Duncan are with Department of Mathematics, University of Kansas, Lawrence, KS 66044, USA, (e-mail: duncan@ku.edu).}
}

\date{First draft: November 2018. This version: April 2019}

\maketitle

\begin{abstract}
This article examines  mean-field-type game problems by means of a direct method.
We provide  various solvable examples beyond the classical linear-quadratic game problems. These  include quadratic-quadratic  games and games with power, logarithmic, sine square, hyperbolic sine square  payoffs. Non-linear state dynamics such as log-state, control-dependent regime switching, quadratic state, cotangent state and hyperbolic cotangent state are considered.  We identify equilibrium strategies and equilibrium payoffs  in state-and-conditional mean-field type feedback form.
It is shown that a simple direct method can be used to solve   broader classes of non-quadratic mean-field-type games under jump-diffusion-regime switching Gauss-Volterra processes which include fractional Brownian motions and multi-fractional Brownian motions. We provide semi-explicit solutions to   the fully cooperative, noncooperative nonzero-sum, and adversarial game problems. 

\end{abstract}


\newpage

\tableofcontents 

\section{Introduction} 
Mean-field-type game theory studies a class of games in which  the payoffs and or state dynamics depend not only on the state-action pairs but also the distribution of them.  
{\color{blue} In mean-field-type games, (i)  a single decision-maker can have a strong impact on the mean-field terms, (ii) the expected payoffs are not necessarily linear with  respect to the state distribution, (iii) the number of decision-makers (``true decision-makers")  is not necessarily infinite.}
 
Games with  non-linearly distribution-dependent quantity-of-interest \cite{alex00,alex01,alex02} are very attractive in terms of applications because 
 the non-linear dependence of the payoff functions in terms of state distribution   allow us to capture risk measures which are functionals of variance, inverse quantiles, and or higher moments. During the past, a significant amount of research on  mean-field-type games has been performed \cite{alex03,alex04,duncan_t1,alex1,alex2, alex3}. In the time-dependent case, the analysis of mean-field-type games  has several  challenges. 
Previous works have devoted   tremendous effort in terms of partial integro-differential  system of equations (PIDEs),  in infinite dimensions, of conditional Liouville, Boltzmann, Kolmogorov or McKean-Vlasov type. At the same time,   an important set of numerical tools have been developed to address the master equilibrium system. However, the current state-of-the-art of numerical schemes is problem-specific and needs to be adjusted properly depending on the underlying problem.  To date, the question of computation  of the master system in the general setting  remains open. This work provides 
explicit solutions of a class of master systems. These explicit solutions can be used to build reference trajectories and  several numerical schemes developed to solve PIDEs  can be tested beyond the linear-quadratic setting.

%\subsection{Direct Method }
%
%\begin{figure}[h!]
%	\begin{center}
%		\resizebox{0.7\textwidth}{!}{
%		\smartdiagram[bubble diagram]{ Direct Method, {\bf 1)} MFTG \\Problem, {\bf 2)} Guess\\ Functional, {\bf 3)} Stochastic \\ Integration \\ , {\bf 4)} Terms \\ Completion, {\bf 5)} Process\\ Identification}
%		}
%	\end{center}
%	\caption{\textcolor{blue}{Direct method and its key   steps. }}
%\label{fig:sum:direct}
%\end{figure}
%
%
%The direct method  consists of five elementary steps. The first step starts by setting the mean-field terms of the problem. The second step consists of the  identification  of a partial guess functional where the coefficient functionals are random and regime switching dependent.  The third step uses the stochastic integration  formula. The fourth step uses  a completion of terms in one-shot optimization for both control actions  and conditional expected value  the control actions of all decision-makers. The fifth and last step uses an algebraic basis of  linearly independent processes to identify the coefficients.  The identification leads  to a (possibly stochastic)  differential system of equations, providing a semi-explicit representation of the solution. These five  elementary steps of the Direct method are displayed in Figure \ref{fig:sum:direct}. 

\subsection{Direct Method for LQ-MFTG}
In the current literature, only  relatively few examples of explicitly solvable mean-field-type game problems are available. The most notable examples are (i)   linear-quadratic mean-field-type games (LQ-MFTG) \cite{duncan_t1}, (ii) linear-exponentiated quadratic mean-field-type games  (LEQ-MFTG) \cite{duncan_risk} ,  (ii) adversarial linear-quadratic mean-field-type games  (minmax LQ, minmax LEQ-MFTG) \cite{duncan_t1}. 
 In LQ-MFTG the base state dynamics has two  components:  drift and noise.
\begin{itemize} \item 
the drift is an affine function of the state, expected value of the state, control action and expected value of the control actions of all decision-makers. The coefficients are regime switching dependent. 
\item the noises are combination of diffusion, Gauss-Volterra, jump, regime-switching process where the noise coefficients  are  affine functions of the state, expected value of the state, control action and expected value of the control actions of all decision-makers. The coefficients are regime switching and jump dependent.
\end{itemize}
To the state dynamics, one can add a common noise which is a diffusion-Gauss-Volterra-jump-regime-switching process. 
The cost functions are polynomial of degree two and include the weighted conditional variances, co-variances between state and control actions of all decision-makers. In addition, the cost functional is not measured perfectly. Only a noisy cost is available.  

This basic model of LQ mean-field-type games captures several interesting features such as heterogeneity, risk-awareness and empathy of the decision-makers. 

To solve LQ-MFTG problems one can use the direct method proposed in Figure \ref{fig:sum:direct}. This solution approach does not require solving the Bellman-Kolmogorov equations or backward-forward stochastic differential equations of Pontryagin's type. 
The proposed direct method can be easily implemented by beginners and engineers who are new to the emerging field of mean-field-type game theory. 

   For this broader  class of  
LQ-MFTG problem one can derive a   semi-explicit  solution under sufficient conditions. The existence of solution to the master system corresponding to the LQ-MFTG problem  can be converted into an existence of solution to a system of ordinary differential equations driven by common noises.
 In some particular cases, these systems are stochastic Riccati systems and extensions of Riccati to include some fractional order terms.

\subsection{Direct Method beyond LQ-MFTG}
 The direct method is not limited to the linear-quadratic case. 
The direct method can be extended  to a class of  LEQ-MFTG, minmax LQ-MFTG  and minmax LEQ-MFTG. In this article, we present several examples to illustrate how the  direct method   addresses non-linear and/or non-quadratic  mean-field-type games. The examples below
go beyond  LQ-MFTG, LEQ-MFTG and minmax LQ problems.

The contributions of this article can be summarized as follows. We provide semi-explicit solution for  classes of mean-field-type game problems presented in Table \ref{table:semi:explicit}. Several  noises are examined: Brownian motion $B$, regime switching $s$, jump process $N$, and Gauss-Volterra process $B_{gv}$. The Gauss-Volterra noise processes are obtained from the integral of a Brownian motion with a suitable kernel function.  In addition, several type of common noises are considered: $s, B_o, N_o,B_{o,gv}.$. We limit ourselves to the class of state-and-conditional mean-field type feedback strategies. The analysis for more general class of strategies is beyond the scope of this article.

{\tiny 
%\begin{table*}. 
\begin{figure}[h!]   \hspace{-3cm}
\textcolor{blue}{  
\begin{tabular}{llll}
\hline 
\textbf{Problem}  &  \textbf{State} & \textbf{Cost}  & \textbf{Noise} \\
\hline
Prop. \ref{log:pro1:semi:explicit} & $\begin{array}{l} \mbox{Drift:~} x \ln x + x \sum_{j \in \mathcal{I}} b_{2j}u_j \\ \mathcal{I}= \{1,\dots,I\} \end{array}$ & $\begin{array}{l} q_i \ln (x) + r_i u_i^{2k}, k>\frac{1}{2} \end{array}$ & $\begin{array}{l} \mbox{Brownian:~} xdB \\ \mbox{Jump:~} \int_{\Theta} x\mu d\tilde{N} \\ \mbox{Switching:~} \tilde{q}_{ss'}  \end{array}$\\ \hline
Prop. \ref{semi:prop2}  & $\begin{array}{l} \mbox{Drift:~}  b_1 x \ln(x) + \sum_{j \in \mathcal{I}} b_{2j} x u_j \\ \mathcal{I}= \{1,\dots,I\} \end{array}$ & $\begin{array}{l} q_i \ln^2(x) + r_iu_i^2 \end{array}$ & $\begin{array}{l} \mbox{Brownian:~} x \sqrt{\ln(x)}\\  \mbox{Switching:~} \tilde{q}_{ss'}  \end{array}$\\ \hline
Prop. \ref{semi:prop3} & $\begin{array}{l} \mbox{Drift:~} \frac{b_1l_1(x)}{l'_1(x)} + \frac{h(x)\sum_{j}b_{2j} u_j}{l'_1(x)} \\ \mathcal{I}= \{1,\dots,I\} \end{array}$ & $\begin{array}{l} q_i l_1(x)+\sum_{j\in \mathcal{I}}r_{ij} l_1(u_j) \end{array}$ & $\begin{array}{l} \mbox{Brownian:~} \sqrt{\frac{ \sigma_1^2+ \sigma_2^2 l_1(x)}{l''_1(x)}} \\ \mbox{Switching:~} \tilde{q}_{ss'}  \end{array}$\\ \hline
Prop. \ref{prop:geo_gauss_volterra} & $\begin{array}{l} \mbox{Drift:~} 2b\sqrt{x} \\
+ x \Big[ \sigma^2 + \sigma_{o}^2 + \sigma_{cogv}^2 +\sigma_{o,cogv}^2   \\
+\int_{\Theta} \mu^2 \nu(d\theta) + \int_{\Theta} \mu_{o}^2 \nu_{o}(d\theta) \Big] \\ \mathcal{I}= \{1,\dots,I\} \end{array}$ & $\begin{array}{l} q_{i}\frac{x^k}{k} + r_{i}\frac{u_i^{2k}}{k} \end{array}$ & $\begin{array}{l} \mbox{Brownian:~} 2x\sigma dB \\ \mbox{Common Noise:~} 2x\sigma_o dB_o \\ \mbox{Jump:~} 2x \int_{\Theta} \mu d\tilde{N} \\ \mbox{Common Jump:~} 2x \int_{\Theta} \mu_o d\tilde{N}_o \\ \mbox{Switching:~} \tilde{q}_{ss'} \\ \mbox{Common G-V:~} 2x\sigma_{o,gv} dB_{o,gv}  \\ \mbox{Gauss-Volterra:~} 2x\sigma_{gv} dB_{gv} \end{array}$\\ \hline
Prop. \ref{semi:prop4} & $\begin{array}{l} \mbox{Drift:~} 0 \\ s(t) \in \mathcal{S}\\ \mathcal{I}= \{1,\dots,I\} \end{array}$ & $\begin{array}{l} r_i (u_i-\bar{u}_i)^2+\bar{r}_i\bar{u}_i^2 + \epsilon_{i}\bar{u}_i \end{array}$ & $\begin{array}{l} \mbox{Switching:~} \tilde{q}_{ss'}(u)\\=\sum_{j}b_{2jss'} (u_j-\bar{u}_j)^2 +
\bar{b}_{2jss'}\bar{u}_j^2 \\+ \sum_{j}{b}_{1jss'}(u_j -\bar{u}_j) + \bar{b}_{1jss'}\bar{u}_j \\+\sum_{j}\bar{b}_{ojss'}, \end{array}$\\ \hline
Prop. \ref{semi:prop5} & $\begin{array}{l} \mbox{Drift:~} \sum_{j \in \mathcal{I}} [q_j (u_j-\bar{u}_j)^2+\bar{q}_j\bar{u}_j^2 \\+ \epsilon_{1j}(u_j - \bar{u}_j) + \bar{\epsilon}_{1j}\bar{u}_j] \\ \mathcal{I}= \{1,\dots,I\} \end{array}$ & $\begin{array}{l} r_i (u_i-\bar{u}_i)^2+\bar{r}_i\bar{u}_i^2 + \bar{\epsilon}_{2i}\bar{u}_i \end{array}$ & $\begin{array}{l} \mbox{Brownian:~} \sigma dB \\ \mbox{Jump:~} \int_{\Theta} \mu d\tilde{N}  \end{array}$\\ \hline
Prop. \ref{semi:prop6} & $\begin{array}{l} \mbox{Drift:~} \frac{1}{2}cot(\frac{x-\bar{x}}{2})+\sum_{j}b_{2j}(u_j-\bar{u}_j)\\
+\frac{1}{2}cot(\frac{\bar{x}}{2})+\sum_{j}\bar{b}_{2j} \bar{u}_j \\ \mathcal{I}= \{1,\dots,I\} \end{array}$ & $\begin{array}{l} ((u_i-\bar{u}_i)^2-q_i)\cos^2(\frac{x-\bar{x}}{4})+q_i  \\
+(\bar{u}_i^2-\bar{q}_i)\cos^2(\frac{\bar{x}}{4})+\bar{q}_i \end{array}$ & $\begin{array}{l} \mbox{Brownian:~} \sigma dB  \\ \mbox{Switching:~} \tilde{q}_{ss'}  \end{array}$\\ \hline
Prop. \ref{semi:prop7} & $\begin{array}{l} \mbox{Drift:~} \frac{1}{2}coth(\frac{x-\bar{x}}{2})+\sum_{j}b_{2j}(u_j-\bar{u}_j)\\ +\frac{1}{2}coth(\frac{\bar{x}}{2})+\sum_{j}\bar{b}_{2j} \bar{u}_j \\ \mathcal{I}= \{1,\dots,I\} \end{array}$ & $\begin{array}{l} ((u_i-\bar{u}_i)^2+q_i)\cosh^2(\frac{x-\bar{x}}{4})-q_i \\ +
 (\bar{u}_i^2+\bar{q}_i)\cosh^2(\frac{\bar{x}}{4})-\bar{q}_i \end{array}$ & $\begin{array}{l} \mbox{Brownian:~} \sigma dB \\ \mbox{Switching:~} \tilde{q}_{ss'}  \end{array}$\\ \hline
Prop. \ref{semi:prop8} & $\begin{array}{l} \mbox{Drift:~} -(u_1-\bar{u}_1)+b_1(x-\bar{x})\\ +b_2\epsilon (x-\bar{x})(u_2-\bar{u}_2)\\ -\bar{u}_1+\bar{b}_{11}\bar{x}+ \bar{b}_{12}\bar{y}+\bar{b}_{13}\bar{z} +\bar{b}_2\bar{u}_2 \bar{x} \\ \mathcal{I}= \{1,\dots,I\} \end{array}$ & $\begin{array}{l} -q \ var(x) -r_1\ var(u_1) \\+ \bar{r}_1  \frac{\bar{u}_1^{\rho}}{\rho} \end{array}$ & $\begin{array}{l} \mbox{Brownian:~} \sigma(x-\bar{x})(u_2-\bar{u}_2)dB \\\mbox{Common noise:~} \bar{\sigma}\bar{x}\bar{u}_2 dB_o \\ \mbox{Switching:~} \tilde{q}_{ss'}   \end{array}$\\ \hline
Prop. \ref{semi:prop9} & $\begin{array}{l} \mbox{Drift:~} b_1 (x-\bar{x})+\sum_{j}b_{2j} (u_j-\bar{u}_j) \\ +\bar{b}_1 \bar{x}+\sum_{j}\bar{b}_{2j} \bar{u}_j \\ \mathcal{I}= \{1,\dots,I\} \end{array}$ & $\begin{array}{l} q_{i}\frac{(x-\bar{x})^{2k_i}}{2k_i} +{r}_i  \frac{(u_i-\bar{u}_i)^{2{k}_i}}{2{k}_i}\\+\bar{q}_{i}\frac{\bar{x}^{2\bar{k}_i}}{2\bar{k}_i} + \bar{r}_i  \frac{\bar{u}_i^{2\bar{k}_i}}{2\bar{k}_i} \end{array}$ & $\begin{array}{l} \mbox{Brownian:~} (x-\bar{x}) \sigma dB \\ \mbox{Jump:~} (x-\bar{x}) \int_{\Theta} \mu d\tilde{N} \\ \mbox{Switching:~}  \tilde{q}_{ss'} \\ \mbox{Gauss-Volterra:~} (x-\bar{x})\sigma_{gv} dB_{gv} \end{array}$ \\  \hline
\end{tabular}
}
%\end{table*}
\caption{\textcolor{blue}{\bf Semi-explicit solutions presented in this article. }}
\label{table:semi:explicit} 
\end{figure}
}

 {\it To  the best of the authors' knowledge this is the first work to provide  semi-explicit solutions  of 
mean-field-type games beyond LQ and under Gauss-Volterra processes.}

\subsection*{Structure}

The rest of the article is structured as follows.  Section \ref{sec:mf:free} presents semi-explicit solutions to some non-linear non-quadratic stochastic differential games. In Section \ref{sec:mf:dependent} we formulate and solve various mean-field-type games with non-quadratic quantity-of-interest and provide semi-explicit solutions using a direct method. Section \ref{sec:mf:gv} presents semi-explicit solutions to some non-quadratic mean-field-type games driven by Gauss-Volterra processes. Numerical examples are presented in Section \ref{sec:mf:num:gv}. The last section summarizes  the work.

\begin{figure}[htb]
\begin{tabular}{|c|c|}
 \hline
Notation & Description  \\ \hline
$B$ & Brownian motion\\ 
$B_o$ & common Brownian motion\\ \hline
$B_{o,gv}$ & Common Gauss-Volterra process\\ 
$B_{gv}$ & Gauss-Volterra process\\ \hline
$\Theta$ & set of jump sizes\\ \hline
$\nu(d\theta),\nu_o(d\theta)$ & Radon measure over $\Theta$ \\ \hline
$\tilde{N}$ &  compensated jump process\\ 
$\tilde{N}_o$ & common compensated jump process\\ \hline
$x$ &   state  \\ \hline
$y$ &   trend  $\int_{-\tau}^0 e^{\lambda t'} x(t+t')dt'$ \\ \hline
$z$ &   delayed state  $x(t-\tau)$ \\ \hline
$\bar{x}$ & conditional state    \\ \hline
$s$ & regime switching process\\ \hline
$\mathcal{I} = \{1,\dots,I\}$ & set of decision-makers\\ \hline
$u_i$ & control action of decision-maker $i$    \\ \hline
$\bar{u}_i$ & conditional  control action  of $i$   \\ \hline 
\end{tabular} \caption{Notations used in this article.}  \label{nota:gv}
\end{figure}

\section*{Preliminary}
We introduce the following notations (see Table \ref{nota:gv}). 
Let $[0,T], \ T>0$  be a fixed time horizon and $(\Omega,\mathcal{F},\mathbb{F}^{B,N, B_{gv},s,B_o,B_{o,gv},N_o}, \mathbb{P})$  be a given filtered
probability space. The filtration $\mathbb{F}=\{{\mathcal{F}}^{B,N,B_{gv},s,B_o,B_{o,gv},N_o}_t,\ 0\leq t \leq T\}$ is the natural filtration of the union of the family $\{B,N, B_{gv},s,B_o,B_{o,gv},N_o\}$ augmented by $\mathbb{P}-$null sets of ${\mathcal F}.$ 
In practice, $B$ is used to capture smaller disturbance, $N$ is used for larger jumps of the system, $B_{gv}$ is used for Gauss-Volterra  processes (including sub- or super diffusion). Let $k\geq 1.$ ${L}^k([0,T]\times \mathcal{S};\mathbb{R})$ is the set of measurable functions $f: \ [0,T]\times \mathcal{S} \rightarrow \mathbb{R}$  such that $\int_0^T |f(t,s)|^k dt<\infty$. 
$\mathcal{L}^k_{\mathbb{F}}([0,T]\times \mathcal{S};\mathbb{R})$ is the set of $\mathbb{F}$-adapted $\mathbb{R}$-valued processes $X(\cdot)$ such that $\mathbb{E} [\int_0^T |X(t)|^k dt ]<\infty.$ The stochastic quantity
$\bar{x}(t)= \mathbb{E} [X(t) | \ \mathcal{F}^{s,B_o,B_{o,gv},N_o}_{t}]$ denotes the conditional expectation of the random variable $X(t)$ with respect to the filtration $\mathcal{F}^{s,B_o,B_{o,gv},N_o}_{t} .$ Note that $\bar{x}$  is a random process. Below, by abuse of notation we use $s(t), x(t)$ for the values $s(t_{-}), x(t_{-})$ inside the jump  processes $N, N_o$ or the regime-switching  process $s$. The set of decision-makers is denoted by $\mathcal{I}=\{ 1,\ldots, I\}.$ An admissible control strategy $u_i$ of the decision-maker $i$  is an $\mathbb{F}$-adapted. We denote the set of all admissible controls
 by $\mathcal{U}_i$:
%
%\begin{align*}
%\begin{array}{l}
%\mathcal{U}_i=\{u_i(\cdot)\in {\mathcal L}^2_{\mathbb{F}}([0,T]\times \mathcal{S};\mathbb{R}); \\
%\,\, u_i(.)\in U_i \,\, a.e.\, t\in [0,T],\,\, \mathbb{P}-a.s. \}.
%\end{array}
%\end{align*}
%
  Decision-maker $i$ chooses a
 control strategy $u_i\in \mathcal{U}_i$ to optimize its performance functional. 
The information structure of the problem under perfect state observation and under common noise observation $(s,B_o,B_{o,gv},N_o).$ 

\subsection{Conditional dynamics of mean-field type}
Consider the following state dynamics of conditional McKean-Vlasov type with time delays, trend, diffusion, jump, regime switching, Gauss-Volterra and common noises.
\begin{equation}
\label{smp00yyt}
\left\{
\begin{array}{lll}
dx= b(t, x, y,z,u,\bar{x},\bar{y}, \bar{z},\bar{u}, m_1,m_2,s) dt  %& \mbox{drift}
\\ + \sigma(t, x, y,z,u,\bar{x},\bar{y}, \bar{z},\bar{u}, m_1,m_2,s) d{B} \\  %& \mbox{diffusion}\\
 + \sigma_{gv}(t, x, y,z,u,\bar{x},\bar{y}, \bar{z},\bar{u}, m_1,m_2,s) d{B}_{gv}\\ % & \mbox{Gauss-Volterra}\\
+\int_{\Theta}\ \mu(t, x, y,z,u,\bar{x},\bar{y}, \bar{z},\bar{u}, m_1,m_2,s,\theta) \tilde{N}(dt,d\theta,s)\\ % & \mbox{jumps}\\ 
\\ + \sigma_{o}(t, x, y,z,u,\bar{x},\bar{y}, \bar{z},\bar{u}, m_1,m_2,s) d{B}_o\\ % & \mbox{common diffusion }\\ 
 + \sigma_{o,gv}(t, x, y,z,u,\bar{x},\bar{y}, \bar{z},\bar{u}, m_1,m_2,s) d{B}_{o,gv}\\ % & \mbox{common Gauss-Volterra }\\
 +\int_{\Theta_o}\ \mu_{o}(t, x, y,z,u,\bar{x},\bar{y}, \bar{z},\bar{u}, m_1,m_2,s,\theta) \tilde{N}_o(dt,d\theta,s)\\ \\ %& \mbox{common jump }\\ \\
x_i(t)=x_{i0}(t),\ \   \ t\in [-\tau_i,0], \   i\in\mathcal{I},\\
% &\mbox{initial state of } \ i, 
% \ i\in \{1,2,\ldots,n \},  \\%& \\  
   s(t)\in\mathcal{S}=\{ 1, 2, \ldots, S\}, \  \\  %%& \mbox{regime switching} \\
   \mathbb{P}(s(t+\epsilon)=s' | s,u)=\int_t^{t+\epsilon} \tilde{q}_{ss'} dt' + o(\epsilon),\ s'\neq s,\ \epsilon>0\\
 \end{array}
\right.
\end{equation}

where \begin{itemize}
\item $\mathcal{I} = \{1,\dots,I\}$ is the set of decision-makers.
 \item $x_i=x_i(t)$ is the basic state at time $t$ of the decision-maker $i$ 
\item  $\tau_{ik}>0$ represents a time delay, 
\item $ y_i=(x_i(t-\tau_{ik}))_{1\leq k\leq K}, $  is a $K-$dimensional delayed state vector,
 \item $ z_i(t)=(\int_{t-\tau_i}^t \lambda_i(dt') \phi_{il}(t,t')x_i(t'))_{l\leq I}$ is the integral state vector of the recent past  state over $[t-\tau_i,t].$
 The trend of the state of decision-maker $i$ is its latest moving averages. $ z_i(t)$  represents the trend of the state of $i$. The process $\phi_{il}(t,t')$ is an $\mathcal{F}_{t'}-$adapted locally bounded process, $\lambda_i$ is a positive and $\sigma-$finite measure on $[-\tau_i, T]$. 
 \item $m_1$ is the distribution states of all the other decision-makers, \item $m_2$ the distribution of actions of all  other decision-makers, 
 \item $x_{i0}$ is a initial deterministic function of  state of $i$ defined on  $[-\tau_i,0].$  
 \item $B$ be a Brownian motion  on $\mathcal{T}=[0,T]$ with suitable dimension. $B_o$ be a  Brownian motion observed by all decision-makers.
 \item $B_{gv}$ be a Gauss-Volterra process on $\mathcal{T}$ with suitable dimension and with integrable kernel $K.$ $B_{o,gv}$ be a  Gauss-Volterra process observed by all decision-makers
 \item ${N}(dt,d\theta,s)$ be a  jump  process with suitable dimension on $\mathcal{T}$ with compensated jump $\tilde{N}(dt,d\theta,s)={N}(dt,d\theta,s)-\nu(d\theta) dt$, $\nu$ is a Radon measure over   $\Theta.$  ${N}_o$ is a   common jump process observed by all decision-makers.
 \item  $s(t)$ is a  regime switching process defined over the finite set  $\mathcal{S}=\{ 1, 2, \ldots, S\}$ with switching rate $\tilde{q}$ satisfying $\tilde{q}_{ss'}>0, s\neq s'$ and $\tilde{q}_{ss}:=-\sum_{s'\neq s}\tilde{q}_{ss'}.$ We use $ \ind_{\{s(t)=s \}}$ to denote the indicator function on the condition $\{s(t)=s \}.$
 \item $ u=(u_i)_{i\in \mathcal{I}}$ is the control strategy profile of all decision-makers. An admissible control strategy $u_i$ of the decision-maker $i$  is an $\mathbb{F}$-adapted process.
 
   \item The  processes $B, B_{gv}, N,  B_o, B_{o,gv},  N_o, s$, are defined in a given filtered probability space $(\Omega,\mathbb{F},\mathcal{F}^{B,N, B_{gv},s,B_o,B_{o,gv},N_o}, \mathbb{P})$
   $ ( \mathbb{F}=\{\mathcal{F}^{B,N, B_{gv},s,B_o,B_{o,gv},N_o}_t\}_{t\in \mathcal{T}}).$ The processes $ B_o, B_{o,gv},  N_o, s$ are common noises assumed to be observable by all decision-makers. All the processes are assumed to be mutually independent.
   \item The coefficient  functionals $b,\sigma,\sigma_{gv}, \mu, \sigma_o,\sigma_{o,gv}, \mu_o$ are of compatible dimensions with $x.$
   
   \item The  quantity $\bar{\omega}(t)= \mathbb{E} [\omega(t) | \ \mathcal{F}^{B_o,B_{o,gv},N_o,s}_{t}]$ denotes the conditional expectation of the random variable $\omega(t)$ with respect to the filtration $\mathcal{F}^{B_o,B_{o,gv},N_o,s}_{t} .$ Note that $\bar{\omega}$  is a random process. We take $\omega\in \{u,x,  y, z\}.$
 By abuse of notation we use $s(t), x(t)$ for the values $s(t_{-}), x(t_{-})$ inside the jump  processes $N, N_o$ or the regime-switching  process $s$. 
 
 \end{itemize}

 Let $f(t,x,s)$ be a twice continuously differentiable function in $x$ and continuously differentiable in time $t$ for each regime $s\in \mathcal{S}.$ Using \cite{prot} and  \cite[Theorem 4.1]{eliott12}, the stochastic integration formula, which is an extended It\^o's formula,  yields 
 \begin{equation}
\label{smp00yyt2}
\left\{
\begin{array}{lll}
f(T,x(T),s(T))\\ 
= f(0,x_0, s_0)+\int_0^T [f_t+ \langle f_x, b\rangle] dt\\
+\frac{1}{2} \int_0^T \langle  f_{xx}\sigma, \sigma  \rangle dt\\
+\frac{1}{2} \int_0^T  \langle  f_{xx}\sigma_{cogv}, \sigma_{cogv}  \rangle dt\\
+\int_0^T \int_{\Theta} [f(x+\mu)-f(x)-\langle \mu, f_x\rangle] \nu(\theta) dt\\
%+\int_0^t \sum_{s'\neq s} [f(., x, s')-f(., x,s)] \tilde{q}_{ss'} dt\\
+\frac{1}{2} \int_0^T  \langle  f_{xx}\sigma_o, \sigma_o  \rangle dt\\
+\frac{1}{2} \int_0^T  \langle  f_{xx}\sigma_{o,cogv}, \sigma_{o,cogv}  \rangle dt\\
+\int_0^T \int_{\Theta} [f(x+\mu_o)-f(x)-\langle \mu_o, f_x \rangle] \nu_o(\theta) dt\\
+\int_0^T \sum_{s'\neq s} [f(., x, s')-f(., x,s)] \tilde{q}_{ss'} dt\\ \\
+  \int_0^T  \langle  f_{x}, \sigma dB \rangle \\
+\int_0^T \langle  f_{x}, \sigma_{gv} dB_{gv}  \rangle \\
+\int_0^T \int_{\Theta} [f(x+\mu)-f(x) \rangle] \tilde{N}(dt,d\theta) \\ \\
+  \int_0^T  \langle  f_{x}, \sigma_o  dB_o \rangle \\
+\int_0^T  \langle  f_{x}, \sigma_{o,gv} dB_{o,gv}  \rangle \\
+\int_0^T \int_{\Theta} [f(x+\mu_o)-f(x) \rangle] \tilde{N}_o(dt,d\theta). \\
 \end{array}
\right.
\end{equation}

Notice that (\ref{smp00yyt2}) applies to  a one-dimensional state as well as to a vector, matrix, tensor, lattice or another object in a Hilbert space. For vectors in an Euclidean space, the inner product is $\langle a,b\rangle=\sum_{l=1}^d a_l b_l,$ for matrices, $\langle a,b\rangle=trace[a^*b] =trace[b^*a],$ where $a^*$ is the transpose of $a.$

\subsection{Direct Method} \label{sec:directm}
\begin{figure}[h!]
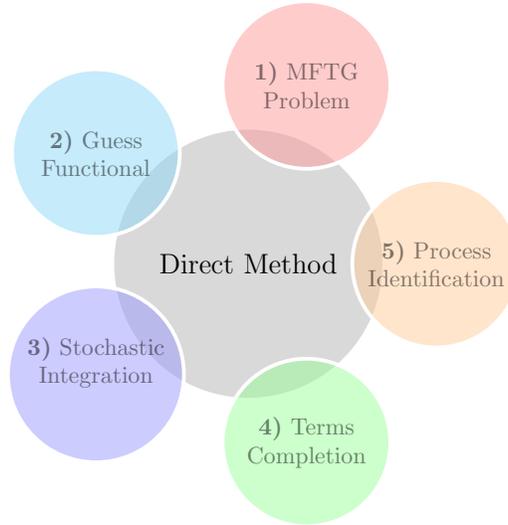

	\begin{center}
		\resizebox{0.6\textwidth}{!}{
		\smartdiagram[bubble diagram]{ Direct Method, {\bf 1)} MFTG \\Problem, {\bf 2)} Guess\\ Functional, {\bf 3)} Stochastic \\ Integration \\ , {\bf 4)} Terms \\ Completion, {\bf 5)} Process\\ Identification}
		}
	\end{center}
	\caption{Direct method and its key   steps.}
\label{fig:sum:direct}
\end{figure}

Consider $I$ decision-makers under perfect state observation $x$ and common noise observation $(B_o, B_{o,gv}, N_o, s).$ 
Given $I$ cost functionals $L_i(x,y,z,u,s)$ associated with (\ref{smp00yyt}), we use  (\ref{smp00yyt2}) in the direct method described as follows. 
The direct method  consists of five elementary steps (see Figure \ref{fig:sum:direct}). 

\begin{itemize}
\item The first step starts by setting the mean-field terms of the problem.
\item The second step consists of the  identification  of a partial guess functional where the coefficient functionals are random and regime switching dependent. For each decision-maker $i$, one needs to identify a guess functional $f_i (t,x,y,z,u,s).$
\item In the third step we compute the difference $L_i-f_i (t,x,y,z,u,s)$ using   the  stochastic integration  formula (\ref{smp00yyt2}).
\item  In the fourth step, we use completion of terms in one-shot optimization for both control actions  and conditional mean-field of  the control actions of all decision-makers. Terms completion  make  $\mathbb{E}[L_i-f_i (t,x,y,z,u,s)]\geq 0$ by matching coefficients. The latter inequality becomes equality iff the optimal control strategies are used. 
\item The fifth and last step uses an algebraic basis of  linearly independent processes to identify the coefficients.  The identification leads  to a (possibly stochastic)  differential system of equations, providing a semi-explicit representation of the solution. The matched coefficients provide  simpler differential systems that are uncoupled with the state.
\end{itemize}

\section{Some Solvable Mean-Field-Free Games} \label{sec:mf:free}
We start with mean-field-free settings where logarithm, logarithm square, Legendre-Fenchel duality, and power payoffs are presented.  The cost functions are not necessarily quadratic and the state dynamics is not necessarily linear.

\subsection{Logarithmic Scale}
Consider a set of decision makers $\mathcal{I} = \{1,\dots,I\} $ interacting in the following non-linear non-quadratic mean-field-free game:
\begin{equation} \label{log:problem1:semi:explicit}
\begin{array}{ll}
\begin{cases}
L_i(x,u) = -q_i(T, s(T))\ln(x(T)) + \int_{0}^{T} \left(-q_i \ln(x) + {\color{blue} r_i u_i^{2k}} \right) dt,\\
\inf_{u_i}~  \mathbb{E} [L_i(x,u)],\\
\text{subject~to}\\
\dx =  \left(b_1x \ln(x) + \sum_{j \in \mathcal{I}} b_{2j} x u_j \right) dt +  {\color{blue} x [ \sigma \dB+\int \mu d\tilde{N}   ]} ,\\
\mathbb{P}(s(t+\epsilon)=s' | s,u)=\int_t^{t+\epsilon} \tilde{q}_{ss'} dt' + o(\epsilon),\ s'\neq s\\
\end{cases}
\end{array}
\end{equation}
and with a given initial condition $x(0) \triangleq x_0>>e, s(0)=s_0\in\mathcal{S},$  $k\geq 1$ is an integer, and 
$\tilde{q}_{ss'}>0, s\neq s'$ and $\tilde{q}_{ss}:=-\sum_{s'\neq s}\tilde{q}_{ss'}.$ 
\eod

\begin{prop} \label{log:pro1:semi:explicit} Assume that $x_0>>e, r_i(.)>\delta>0, q_i(.)\geq 0, \ \mu(\theta)\geq 0, $ 
$\int_{\Theta} [\ln(1+\mu(\theta))-\mu(\theta)]\nu(d\theta) <\infty.$
The non-linear non-quadratic mean-field-free Nash equilibrium and the corresponding equilibrium cost are given by:
\begin{align*}
u_i^* &=\sum_{s\in\mathcal{S}} \ind_{\{s(t)=s \}}  [\frac{1}{2k}\frac{\alpha_i b_{2i}}{r_i}]^{\frac{1}{2k-1}},\\
%%u_i^* &= -\frac{1}{2}\frac{\alpha_i b_{2i}}{r_i},\\
\mathbb{E}[L_i(x,u^*)] &= \mathbb{E}[-\alpha_i(0, s_0) \ln(x_0) + \delta_i(0, s_0)],
\end{align*}
where $\alpha_i$ and $\delta_i$ satisfies the following differential equations:
\begin{equation} \label{fieq00}\begin{array}{ll}
\dot{\alpha}_i + q_i + \alpha_ib_1 +\sum_{s'}[{\alpha}_i(t,s')-{\alpha}_i(t,s)]\tilde{q}_{ss'}=0,\\
\dot{\delta}_i + \alpha_i[\frac{\sigma^2}{2}-\int_{\theta\in \Theta}[\ln(1+\mu(\theta))-\mu(\theta)]\nu(d\theta)]\\
+ \sum_{s'}[{\delta}_i(t,s')-{\delta}_i(t,s)]\tilde{q}_{ss'}
-(2k-1) r_i( \frac{1}{2k r_i}b_{2i}\alpha_i)^{\frac{2k}{2k-1}} \\ -\alpha_i\sum_{j\neq i} b_{2j} [\frac{1}{2k}   \frac{\alpha_jb_{2j}}{r_j}   ]^{\frac{1}{2k-1}}
=0, \end{array}
\end{equation}
where $\alpha_i(T,s) = q_i(T,s)$, and $\delta_i(T,s) = 0$\eod. 
\end{prop}. 

\textbf{Proof.} Consider the following guess functional:
\begin{align*}
f_i(t,x,s) = -\alpha_i \ln(x) + \delta_i.
\end{align*}

%It follows that
%\begin{align*}
%\partial_t f_i(t,x) &= \dot{\alpha}_i \ln(x) + \dot{\delta}_i,\\
%\partial_x f_i(t,x) &= \frac{\alpha_i}{x},\\
%\partial_{xx} f_i(t,x) &= -\frac{\alpha_i}{x^2}.
%\end{align*}
By applying  It\^o's formula for jump-diffusion-regime switching processes, the gap between the cost and the guess functional $\mathbb{E}[L_i(x,u)- f_i(0,x_0)]$ can be computed and it is given by

\begin{equation}
\begin{array}{ll}
\mathbb{E}[L_i(x,u)- f_i(0,x_0,s_0)] \\ 
= \mathbb{E} \left(-q_i(T,s(T))+\alpha_i(T, s(T))\right) \ln(x(T)) + \delta_i(T, s(T))\\
 \mathbb{E} \int_{0}^{T}
-\{ \dot{\alpha}_i + q_i + \alpha_ib_1 +\sum_{s'}[{\alpha}_i(t,s')-{\alpha}_i(t,s)]\tilde{q}_{ss'}\} ln(x) dt\\
+\int_{0}^{T}\dot{\delta}_i +\frac{\sigma^2}{2}\alpha_i+ \sum_{s'}[{\delta}_i(t,s')-{\delta}_i(t,s)]\tilde{q}_{ss'} dt\\
+\int_{0}^{T} -(2k-1) r_i( \frac{1}{2k r_i}b_{2i}\alpha_i)^{\frac{2k}{2k-1}} 
-\alpha_i\sum_{j\neq i} b_{2j} [\frac{1}{2k}   \frac{\alpha_jb_{2j}}{r_j}   ]^{\frac{1}{2k-1}} \\
-\int_{0}^{T}\alpha_i\int_{\theta\in \Theta}[\ln(1+\mu(\theta))-\mu(\theta)]\nu(d\theta) dt\\
+\mathbb{E} \int_{0}^{T}[-b_{2i}\alpha_i  u_i+ r_iu_i^{2k}  + (2k-1) r_i( \frac{1}{2k r_i}b_{2i}\alpha_i)^{\frac{2k}{2k-1}} 
 ] dt,
\end{array}
\end{equation}
Noting that $$[-b_{2i}\alpha_i  u_i+ r_iu_i^{2k} + (2k-1) r_i( \frac{1}{2k r_i}b_{2i}\alpha_i)^{\frac{2k}{2k-1}}  ] \geq 0$$ with equality iff 
$u_i=u_i^*:=[\frac{1}{2k}\frac{\alpha_i b_{2i}}{r_i}]^{\frac{1}{2k-1}},$  the announced result follows. \eod

{\color{blue} Notice that the differential system (\ref{fieq00}) has a unique solution:  the system in $\alpha$ is linear and the system in $\delta$ is obtained by integration.   $\ln(x)$ is well-defined because the state $x$ stays positive in $[0,T]$ almost surely  if one starts at $x_0 >> e.$ }

\begin{remark}
For $k=1$ the system reduces to the following ordinary differential equations:

\begin{equation}\begin{array}{ll}
u_i^* = [\frac{1}{2}\frac{\alpha_i b_{2i}}{r_i}],\\

\mathbb{E}[L_i(x,u^*)] = \mathbb{E}[\alpha_i(0) \ln(x_0)] + \delta_i(0)],\\
\dot{\alpha}_i + q_i + \alpha_ib_1 +\sum_{s'}[{\alpha}_i(t,s')-{\alpha}_i(t,s)]\tilde{q}_{ss'}=0,\\
\dot{\delta}_i +\frac{\sigma^2}{2}\alpha_i+ \sum_{s'}[{\delta}_i(t,s')-{\delta}_i(t,s)]\tilde{q}_{ss'}\\
- \frac{1}{4 r_i} b_{2i}^2\alpha_i^{2}
-\frac{1}{2} \alpha_i\sum_{j\neq i} b^2_{2j} \frac{\alpha_j}{r_j}   \\
-\alpha_i\int_{\theta\in \Theta}[\ln(1+\mu(\theta))-\mu(\theta)]\nu(d\theta)=0
\end{array}
\end{equation}

\end{remark}
\subsection{Logarithm square}
Consider the following non-linear non-quadratic mean-field-free game:
\begin{equation} \label{semi:problem2}
 \begin{array}{ll} 
\begin{cases}
L_i(x,u) = q_i(T,s(T))\ln^2(x(T)) + \int_{0}^{T} \left(q_i \ln^2(x) + r_iu_i^2 \right) dt,\\
\inf_{u_i}~ \mathbb{E}[L_i(x,u)],\\
\text{subject~to}\\
\dx =  \left( b_1 x \ln(x) + \sum_{j \in \mathcal{I}} b_{2j} x u_j \right) dt  ,\\
x(0) \triangleq x_0>>e,\\
\mathbb{P}(s(t+\epsilon)=s' | s,u)=\int_t^{t+\epsilon} \tilde{q}_{ss'} dt' + o(\epsilon),\ s'\neq s\\
\end{cases}\end{array}
\end{equation}
 \eod

\begin{prop} \label{semi:prop2}
 Assume that  $x_0>>e,\  q_i(t,s)\geq 0,  r_i(t,s)>\delta>0,$ and $\int_{\theta\in \Theta}[\ln(1+\mu(\theta))-\mu(\theta)]\nu(d\theta)<\infty.$
The non-linear non-quadratic mean-field-free Nash equilibrium and corresponding optimal cost are given by:  
\begin{align*}
u_i^* &= - \sum_{s\in\mathcal{S}} \ind_{\{s(t)=s \}}  \frac{\alpha_ib_{2i}}{r_i} \ln(x),\\
\mathbb{E}[L_i(x,u^*)] &= \mathbb{E}[\alpha_i(0) \ln^2(x_0)],
\end{align*}
where $\alpha_i$ satisfies the following differential equation:
\begin{align*}
\dot{\alpha}_i = -q_i   -2 b_1 \alpha_i + 2\alpha_i \sum_{j \in \mathcal{I}\setminus\{i\}} \alpha_j \frac{b_{2j}^2}{r_j} \\
+ \alpha_i^2 \frac{b_{2i}^2}{r_i}+\sum_{s'\in \mathcal{S}} [\alpha_i(t,s')-\alpha_i(t,s)]\tilde{q}_{ss'},\\
 \alpha_i(T,s) = q_i(T,s),\\ 
\end{align*}
 \eod
\end{prop}

\textbf{Proof.} Consider the following guess functional:
\begin{align*}
f_i(t,x) = \alpha_i \ln^2(x).
\end{align*}
%where
%\begin{align*}
%\partial_t f_i(t,x) &= \dot{\alpha}\ln^2(x),\\
%\partial_x f_i(t,x) &= 2\alpha_i \frac{\ln(x)}{x},\\
%\partial_{xx} f_i(t,x) &= 2\alpha_i \frac{\left(1-\ln(x)\right)}{x^2}.
%\end{align*}
Applying the It\^o's formula yields
\begin{align*}
&f_i(T,x(T)) - f_i(0,x_0) = \int_{0}^{T} \dot{\alpha}_i\ln^2(x) dt \\
&+ \int_{0}^{T} 2\alpha_i \ln(x) \left( b_1 \ln(x) + \sum_{j \in \mathcal{I}} b_{2j} u_j \right) dt \\
&+\int_{0}^{T}\sum_{s'\in \mathcal{S}} [\alpha_i(t,s')-\alpha_i(t,s)]\tilde{q}_{ss'} \ln^2(x) dt
%&+\int_{0}^{T} \alpha_i \left(1-\ln(x)\right) \ln(x) dt + \int_{0}^{T} 2\alpha_i \ln(x) \sqrt{\ln(x)} \dB.
\end{align*}
Thus, the gap $\mathbb{E}[L_i(x,u)- f_i(0,x_0)]$ is given by
\begin{align*}
&\mathbb{E}[L_i(x,u)- f_i(0,x_0)] = \mathbb{E}\left(q_i(T)-\alpha_i(T)\right) \ln^2(x(T))\\
&+ \mathbb{E} \int_{0}^{T} q_i \ln^2(x) dt + \mathbb{E} \int_{0}^{T} \dot{\alpha}_i \ln^2(x)  dt \\
&+ \mathbb{E} \int_{0}^{T} \left(2\alpha_i b_1 \ln^2(x) + 2\alpha_i \ln(x) \sum_{j \in \mathcal{I}\setminus\{i\}} b_{2j} u_j \right) dt \\
&+\mathbb{E} \int_{0}^{T} r_i\left(u_i^2 + 2\alpha_i \frac{\ln(x) b_{2i}}{r_i} u_i\right) dt\\
&+\int_{0}^{T}\sum_{s'\in \mathcal{S}} [\alpha_i(t,s')-\alpha_i(t,s)]\tilde{q}_{ss'} \ln^2(x) dt
\end{align*}
By performing square completion one obtains
\begin{align*}
\left(u_i+ \alpha_i \frac{\ln(x) b_{2i}}{r_i}\right)^2 - \alpha_i^2 \frac{\ln^2(x) b_{2i}^2}{r_i^2} = u_i^2 + 2\alpha_i \frac{\ln(x) b_{2i}}{r_i}u_i,
\end{align*}
then, 
\begin{align*}
&\mathbb{E}[L_i(x,u)- f_i(0,x_0)] = \left(q_i(T)-\alpha_i(T)\right) \ln^2(x(T)) \\
&+ \mathbb{E} \int_{0}^{T} q_i \ln^2(x) dt + \mathbb{E} \int_{0}^{T} \dot{\alpha}_i\ln^2(x) dt \\
&+ \mathbb{E} \int_{0}^{T} \left(2\alpha_i b_1 \ln^2(x) - 2\alpha_i \ln^2(x) \sum_{j \in \mathcal{I}\setminus\{i\}} \alpha_j \frac{b_{2j}^2}{r_j} \right) \ dt \\
& + \mathbb{E} \int_{0}^{T} r_i\left(u_i+ \alpha_i \frac{\ln(x) b_{2i}}{r_i}\right)^2 dt \\
&- \mathbb{E} \int_{0}^{T} \alpha_i^2 \frac{\ln^2(x) b_{2i}^2}{r_i} dt\\
&+\mathbb{E} \int_{0}^{T}\sum_{s'\in \mathcal{S}} [\alpha_i(t,s')-\alpha_i(t,s)]\tilde{q}_{ss'} \ln^2(x) dt
\end{align*}
Finally, the announced result is obtained by minimizing the terms. \eod

\subsection{Legendre-Fenchel}
We  consider a convex running loss functions $l_1,l_2.$  

\begin{equation} \label{problem:fenchel}
\begin{cases}
L_i = q_{iT}l_1(x_T)+\int_0^T q_i l_1(x)+\sum_{j\in \mathcal{I}}r_{ij} l_2(u_j) dt,\\
\inf_{u_i} \mathbb{E}[L_i], \\
\mbox{subject to }\\
\mathbb{P}(s(t+\epsilon)=s' | s)=\int_t^{t+\epsilon} \tilde{q}_{ss'} dt' + o(\epsilon),\ s'\neq s\\
dx=[ \frac{b_1l_1(x)}{l'_1(x)} + \frac{h(x)\sum_{j}b_{2j} u_j}{l'_1(x)}]dt+\sqrt{\frac{ \sigma_1^2+ \sigma_2^2 l_1(x)}{l''_1(x)}}dB,\\
x(0)=x_0, \ s(0)=s_0
\end{cases}
\end{equation}
where $\tilde{q}_{ss'}>0, s\neq s'$ and $\tilde{q}_{ss}:=-\sum_{s'\neq s}\tilde{q}_{ss'}.$ 
Recall that the Legendre-Fenchel transform of $l,$ is given by  $$-l^*(x)=\inf_{u}\{ l(u)-xu\}.$$

\begin{prop}  \label{semi:prop3}
Assume that   $ l_1, l_2, l''_2, r_{ij}, q_i \ $ are  positive.  Then, the game problem (\ref{problem:fenchel}) has a solution:

\begin{equation}
\begin{array}{ll}
u^*_i= \sum_{s\in\mathcal{S}} \ind_{\{s(t)=s \}} (l^*_2)'(-\frac{b_{2i}\alpha_i}{r_{ii}} h(x)), \\
\mathbb{E}[L_i(x,u)] =\mathbb{E}[ \alpha_i(0,s_0) l(x_0)+\delta_i(0,s_0)],\\
\end{array}
\end{equation}
with 
\begin{equation}
\begin{array}{ll}
 \dot{\alpha}_i +q_i+\alpha_i (b_1+\frac{\sigma_2^2}{2})+\sum_{s'}(\alpha_i(t,s')-\alpha_i(t,s))\tilde{q}_{ss'}\\
 -\eta_{ii} + \sum_{j\neq i} \eta_{ij} +\alpha_ib_{2j}\gamma_j 
 =0 ,\\
\dot{\delta}_i+\frac{\sigma_1^2}{2}+\sum_{s'}(\delta_i(t,s')-\delta_i(t,s))\tilde{q}_{ss'}=0\\
\end{array}
\end{equation}
where  
\begin{equation}
\begin{array}{ll}
r_{ii} l^{*}_2(-\frac{b_{2i}\alpha_i}{r_{ii}} h(x)) =\eta_{ii} l_1(x),\\
 r_{ij} l_2(u^*_j)= \eta_{ij} l_1(x),\\
h(x)(l_2^*)'[ -\frac{\alpha_j b_{2j}}{r_{jj}} h(x)]= \gamma_j l_1(x),\\
\end{array}
\end{equation}
These conditions are fulfilled by choosing for example $$l_2(y)= y h(y),\ l_1=\kappa l_2, h(x)=\frac{x^{2k-1}}{2k}.$$ \eod
\end{prop}

{\bf Proof}
Step 1: we observe that the structure of the problem is mainly driven by the evolution of the function $l_1.$ 

Step 2: Inspired the nature of the problem, we propose a guess functional in the form of $l_1$ with deterministic coefficients $\alpha_i,\delta_i.$
Let  $f_i(t,x,s)=\alpha_i(t,s) l_1(x)+\delta_i(t,s).$

Step 3: We apply stochastic integration formula for diffusion-regime switching to obtain the difference between the cost and the guess functional as
\begin{equation}
\begin{array}{ll}

\mathbb{E}[L_i(x,s,u)- f_i(0,x_0,s_0)]\\
= \mathbb{E} (q_i(T,s(T))-\alpha_i(T,s(T))) l_1(x(T))+(0-\delta_i(T,s(T)))\\
+ \mathbb{E} \int_0^T  [\dot{\alpha}_i +q_i+\alpha_i (b_1+\frac{\sigma_2^2}{2})\\
+\sum_{s'}(\alpha_i(t,s')-\alpha_i(t,s))\tilde{q}_{ss'}]l(x)\\
+\dot{\delta}_i+\frac{\sigma_1^2}{2}+\sum_{s'}(\delta_i(t,s')-\delta_i(t,s))\tilde{q}_{ss'}\\
 - r_{ii} l^{*}_2(-\frac{b_{2i}\alpha_i}{r_{ii}} h(x))   +\sum_{j\neq i}  \{ r_{ij} l_2(u^*_j) +\alpha_ib_{2j}h(x) u^*_j\}\\
+  \{ r_{ii}[ l_2(u_i) +\frac{\alpha_i b_{2i}h(x)}{r_{ii}} u_i]  + r_{ii} l_2^{*}(-\frac{\alpha_i b_{2i}h(x)}{r_{ii}} )\} dt\\
\end{array}
\end{equation}

Step 4: Observing that 
  $$\{  [l_2(u_i) +\frac{\alpha_i b_{2i}h(x)}{r_{ii}} u_i]  + l_2^{*}(-\frac{\alpha_i b_{2i}h(x)}{r_{ii}} )\}\geq 0,$$ with equality iff $u_i=u_i^*,$ 
the  one-shot optimization provides $u_j^*=(l_2^*)'[- \frac{\alpha_j b_{2j}}{r_{jj}} h(x)],$  

Step 5: By identification of processes, the announced result follows. This completes the proof. 
\eod

\subsection{Geometric Gauss-Volterra Game}

The Gauss-Volterra processes are singular integrals of a standard Brownian motion and 
include   (i) fractional Brownian motions, (ii) Liouville fractional Brownian motions,  and (iii) multi-fractional Brownian motions. 
The difficulty of finding semi-explicit solution is significantly increased if the noise process and thereby the state process is driven by non-Markov processes or non-martingales.
Let $B_{gv}$ be a Gauss-Volterra process with  zero mean and covariance 
$$c(t,t')=\mathbb{E} [B_{gv}(t)B_{gv}(t')]=\int_0^{\min(t,t')} K(t,t'')K(t',t'') dt''.$$  
%The kernel $K$ has the following properties:
%\begin{itemize}
%\item $K(0,0)=0, $  $K(t,t')=0$ if $0<t<t'.$ for every time $t,$ $\int_0^t K^2(t,t')dt' <+\infty.$ For each $T'>0$  there are positive constants  $c_1,c_2$ such that 
%$\int_0^{T'} [K(t,t'')-K(t',t'')]^2 dt''\leq c_1|t-t'|^{c_2},\  (t,t')\in [0,T']^2.$  $t\mapsto K(t,t')$ is differentiable in the first variable in $0<t'<t<+\infty.$ Both $K, K_t$ are continuous and $K(t'_+,t')=0$  for each $t'\in [0,+\infty).$
%\item %$ |K_t(t,t')|\leq c_T(t-t')^{c_3-1}(\frac{t}{t'})^{\alpha}$ and
%  $\int_0^t K^2(t,t')dt'\leq c_T(t-s)^{1-2c_3},$ on the set $0<s<t<T$ for some constants $c_T>0$ and $c_3\in (0,\frac{1}{2}).$
%\end{itemize}
%The kernel $K$ has causality and continuity properties from the above  conditions. 
%
The kernel $K$ is assumed to have causality, continuity and integrability properties as in \cite{refpp16}. 
The variance of the process $\int_0^t \sigma_{gv}(t',s) dB_{gv}(t')$ is given by

 $$\sigma_{cogv}^2=\frac{d}{dt}\bigg[\int_0^t  \Bigg\{ K(t'_+,t')\sigma_{gv}(t')+\int_{t'}^{t}\sigma_{gv}(t'')K(t',t'')dt''\Bigg\}^2  dt'\bigg].$$

Consider the following geometric Gauss-Volterra game with unobserved processes $B_o,B_{o,gv},N_o$ which are assumed to be independent.
\begin{align}
\label{eq:geo_volterra}
\left\{
\begin{array}{l}
\inf_{u_i} ~\mathbb{E} L_i(x,u) = \mathbb{E} q_{iT}\frac{x_T^k}{k} + \int_{0}^{T} \left( q_{i}\frac{x^k}{k} + r_{i}\frac{u_i^{2k}}{k} \right) dt,\\
\mathrm{subject~to}\\
dx = \Big( 2b\sqrt{x} + x \Big[ \sigma^2 + \sigma_{o}^2 + \sigma_{cogv}^2 +\sigma_{o,cogv}^2   \\
+\int_{\Theta} \mu^2 \nu(d\theta) + \int_{\Theta} \mu_{o}^2 \nu_{o}(d\theta) \Big] \Big) dt\\
+ 2 x \Big[  \sigma dB + \sigma_{o} dB_{o} +\int_{\Theta} \mu d\tilde{N} + \int_{\Theta} \mu_{o} d\tilde{N}_o  \\
+ \sigma_{gv} dB_{gv} + \sigma_{o,gv} dB_{o,gv} \Big],
\end{array}
\right.
\end{align}
where $b = b_1 \sqrt{x} + \sum_{j\in \mathcal{I}}{b_{2j}} u_j,$ and $\sigma$, $\sigma_{o}$, $\sigma_{gv}$, $\sigma_{o,gv}$, $\mu$, and $\mu_{o}$ are real valued and  regime-switching dependent functions $s(t)$ with $\tilde{q}_{ss'}>0,\  s\neq s'$ and $\tilde{q}_{ss}:=-\sum_{s'\neq s}\tilde{q}_{ss'}.$ 

\begin{prop}
	\label{prop:geo_gauss_volterra} Assume that   $q_{i}\geq 0,\  r_{i}>\delta>0,\  x_0>0, k\geq 1.$
	The mean-field-free equilibrium for the Geometric Gauss-Volterra Game in \eqref{eq:geo_volterra} is given by
	\begin{align*}
	%u_i^* = \left(- \frac{\alpha_i b_{2i}  y^{\left(k-\frac{1}{2}\right)}}{r_i}\right)^{\frac{2k}{2k-1}},
	u_i^* = \sum_{s\in\mathcal{S}} \ind_{\{s(t)=s \}}  \left(- \frac{\alpha_i b_{2i} }{r_i}\right)^{\frac{1}{2k-1}} \sqrt{x},\\
	\mathbb{E} L_i(x,u^*) =\mathbb{E} \alpha_i(0,s_0) \frac{x_0^k}{k},
	\end{align*}
	where $\alpha_i$ satisfies the following differential equation:
	\begin{equation} \label{fieq04}
	\begin{array}{l}
	 \dot{\alpha}_i +q_i + k \alpha_i ( [ \sigma^2 + \sigma_{o}^2 + \sigma_{cogv}^2 +\sigma_{o,cogv}^2 ]\\
	+\int_{\Theta} \mu^2 \nu(d\theta) + \int_{\Theta} \mu_{o}^2 \nu_{o}(d\theta))\\
	+ 2kb_1 \alpha_i + 2k(k-1) \alpha_i ( \sigma^2 + \sigma_{o}^2 + \sigma_{cogv}^2 +\sigma_{o,cogv}^2 \\
	+\int_{\Theta} \mu^2 \nu(d\theta) + \int_{\Theta} \mu_{o}^2 \nu_{o}(d\theta))\\
	- \left(- \frac{\alpha_i b_{2i} }{r_i}\right)^{\frac{2k}{2k-1}} r_i (k-\frac{1}{2}) \\
	+ 2 \alpha_i k  \sum_{j \ne i} b_{2j}  \left(- \frac{\alpha_j b_{2j} }{r_j}\right)^{\frac{1}{2k-1}} \\
	+ \sum_{s'}(\alpha_i(t,s')-\alpha_i(t,s))\tilde{q}_{ss'} =0,
	\end{array}
	\end{equation}
	\eod
\end{prop}  

\textbf{Proof.}  
This proof is developed following a direct method. 

\textit{Step 1:} Observe that the problem is a  mean-field free problem driven by $\sqrt{x}$.

\textit{Step 2:} Based on the structure of the problem we propose the following guess functional: 
\begin{align*}
f_i(t,x,s) = \alpha_i(t,s) \frac{x^k}{k},
\end{align*} 
\textit{Step 3:} We apply stochastic integration formula for jump-diffusion-regime-switching Gauss-Volterra and common noises to compute the difference between the costs and the guess functionals, i.e.,
\begin{align*}
\begin{array}{l}
\mathbb{E}[L_i(x,u) - f_i(0,s_0)] = \mathbb{E} (q_{i,T}-\alpha_i(T,s(T)))\frac{x^k(T)}{k} \\
+ \mathbb{E}\int_{0}^{T} ( q_i + \dot{\alpha}_i) \frac{x^k}{k} + k \alpha_i \frac{x^{k}}{k} ( [ \sigma^2 + \sigma_{o}^2 + \sigma_{cogv}^2 +\sigma_{o,cogv}^2 ]\\
+ \int_{\Theta} \mu^2 d\nu(\theta) + \int_{\Theta} \mu_{o}^2 d\nu_{o}(\theta))\\
+ 2kb_1 \alpha_i \frac{x^{k}}{k} \\
+ 2k(k-1) \alpha_i \frac{x^{k}}{k} ( \sigma^2 + \sigma_{o}^2 + \sigma_{cogv}^2 +\sigma_{o,cogv}^2 \\
+\int_{\Theta} \mu^2 d\nu(\theta) + \int_{\Theta} \mu_{o}^2 d\nu_{o}(\theta))\\
+ [ r_i\frac{u_i^{2k}}{k} + 2 \alpha_i x^{k-\frac{1}{2}} b_{2i} u_i ]
+ 2 \alpha_i x^{k-\frac{1}{2}}  \sum_{j \ne i} b_{2j} u_j\\
+ \frac{x^k}{k} \sum_{s'}(\alpha_i(t,s')-\alpha_i(t,s))\tilde{q}_{ss'},
\end{array}
\end{align*}

\textit{Step 4:} we perform terms completion:
\begin{align*}
%u_i^* = \left(- \frac{\alpha_i b_{2i}  y^{\left(k-\frac{1}{2}\right)}}{r_i}\right)^{\frac{2k}{2k-1}},
u_i^* =  \left(- \frac{\alpha_i b_{2i} }{r_i}\right)^{\frac{1}{2k-1}} \sqrt{x},
\end{align*}

\textit{Step 5:} We perform process identification after having replaced back the optimal control inputs in the gap $\mathbb{E}[L_i(x,u) - f_i(0,s_0)]$, i.e.,
\begin{align*}
\begin{array}{l}
\mathbb{E}[L_i(x,u) - f_i(0,s_0)] = (q_{i,T}-\alpha_i(T,s(T)))\frac{x^k(T)}{k} \\
+ \int_{0}^{T} ( q_i + \dot{\alpha}_i) \frac{x^k}{k} + k \alpha_i \frac{x^{k}}{k} ( [ \sigma^2 + \sigma_{o}^2 + \sigma_{cogv}^2 +\sigma_{o,cogv}^2 ]\\
+\int_{\Theta} \mu^2 d\nu(\theta) + \int_{\Theta} \mu_{o}^2 d\nu_{o}(\theta))\\
+ 2kb_1 \alpha_i \frac{x^{k}}{k} \\
+ 2k(k-1) \alpha_i \frac{x^{k}}{k} ( \sigma^2 + \sigma_{o}^2 + \sigma_{cogv}^2 +\sigma_{o,cogv}^2 \\
+\int_{\Theta} \mu^2 d\nu(\theta) + \int_{\Theta} \mu_{o}^2 d\nu_{o}(\theta))\\
+ [ r_i\frac{u_i^{2k}}{k} + 2 \alpha_i x^{k-\frac{1}{2}} b_{2i} u_i + \left(- \frac{\alpha_i b_{2i} }{r_i}\right)^{\frac{2k}{2k-1}} r_i (1-\frac{1}{2k})  x^k]\\
+ \Big[ -k \left(- \frac{\alpha_i b_{2i} }{r_i}\right)^{\frac{2k}{2k-1}} r_i (k-\frac{1}{2}) \\
+ 2 \alpha_i k  \sum_{j \ne i} b_{2j}  \left(- \frac{\alpha_j b_{2j} }{r_j}\right)^{\frac{1}{2k-1}} \Big] \frac{x^k}{k}\\
+ \frac{x^k}{k} \sum_{s'}(\alpha_i(t,s')-\alpha_i(t,s))\tilde{q}_{ss'} dt,
\end{array}
\end{align*}
Finally, the announced result is obtained by minimizing terms, which completes the proof. \eod

{\color{blue} Notice that under the conditions: $q_i\geq 0,  r_{i}>\delta>0,   i \in \mathcal{I}, k\geq 1$ the differential system (\ref{fieq04}) has a  positive solution.  }

When $s, B_o,B_{o,gv},N_o$ are  noises observed by all decision-makers (observed common noises), the ordinary differential system in $\alpha$ becomes a stochastic differential system driven by the union of events with $B_o,B_{o,gv},N_o.$
 
\begin{equation}
	\begin{array}{l}
	 -d{\alpha}_i 
	 =q_i  dt+ k \alpha_i (  \sigma^2 + \sigma_{o}^2 + \sigma_{gv}^2 +\sigma_{o,gv}^2 \\
	+\int_{\Theta} \mu^2 \nu(d\theta) + \int_{\Theta} \mu_{o}^2 \nu_{o}(d\theta)) dt\\
	+ 2kb_1 \alpha_i + 2k(k-1) \alpha_i ( \sigma^2 + \sigma_{o}^2 + \sigma_{gv}^2 +\sigma_{o,gv}^2  dt\\
	+\int_{\Theta} \mu^2 \nu(d\theta) + \int_{\Theta} \mu_{o}^2 \nu_{o}(d\theta)) dt\\
	- \left(- \frac{\alpha_i b_{2i} }{r_i}\right)^{\frac{2k}{2k-1}} r_i (k-\frac{1}{2}) dt \\
	+ 2 \alpha_i k  \sum_{j \ne i} b_{2j}  \left(- \frac{\alpha_j b_{2j} }{r_j}\right)^{\frac{1}{2k-1}} dt\\
	+ \sum_{s'}(\alpha_i(t,s')-\alpha_i(t,s))\tilde{q}_{ss'}  dt\\
	%%cross variation
	{\color{blue}
	+ 2 k \Big[  {\alpha}_{i,B_o} \sigma_{o}  + \int_{\Theta} \mu_{o}  {\alpha}_{i,N_o} \nu_{o}(d\theta)
 +cov[\sigma_{o,gv}dB_{gv}, {\alpha}_{i,B_{gv}} dB_{gv} ] \Big]dt}\\
	{\color{blue} + 2 k \alpha_i\Big[  \sigma dB + \sigma_{o} dB_{o} +\int_{\Theta} \mu d\tilde{N} + \int_{\Theta} \mu_{o} d\tilde{N}_o}  \\ 
	{\color{blue}+ \sigma_{gv} dB_{gv} + \sigma_{o,gv} dB_{o,gv} \Big]}
	\end{array}
	\end{equation}
	with the terminal condition $q_i(T,s)\geq 0$ being $\mathcal{F}^{B_o,B_{o,gv},N_o}$-measurable random coefficient.

\section{Some Solvable Mean-Field-Type Games}  \label{sec:mf:dependent}
\subsection{Control-dependent switching MFTG }
In most continuous time MFTG models with regime switching considered in the literature it is assumed that the switching rate $\tilde{q}_{ss'}$ is control-independent. In this subsection, we provide an example with control-dependent switching rate $\tilde{q}_{ss'}(u)$ in which the MFTG problem can be solved  semi-explicitly. 
\begin{equation}
\begin{cases}
L_i = q_{iT}+\int_0^T r_i (u_i-\bar{u}_i)^2+\bar{r}_i\bar{u}_i^2 + \epsilon_{i}\bar{u}_i dt,\\
\inf_{u_i} \mathbb{E}[L_i], \\
\mbox{subject to }\\
\mathbb{P}(s(t+\epsilon)=s' | s,u)=\int_t^{t+\epsilon} \tilde{q}_{ss'}(u) dt' + o(\epsilon),\ s'\neq s\\
\tilde{q}_{ss'}(u)=\sum_{j}b_{2jss'} (u_j-\bar{u}_j)^2+
\bar{b}_{2jss'}\bar{u}_j^2\\ + \sum_{j}{b}_{1jss'}(u_j -\bar{u}_j) + \bar{b}_{1jss'}\bar{u}_j+\bar{b}_{ojss'},\\

\end{cases}
\end{equation}
where $b_{kjss'}>0$ for $s'\neq s$ and $\mathbb{E}{b}_{1jss'}=0.$

{\color{black}
\begin{prop} \label{semi:prop4}  Assume that $ r_i, \bar{r}_i>0.$
The equilibrium strategy is
\begin{equation}
\begin{array}{ll}
u_i^*=  %\sum_{s\in\mathcal{S}} \ind_{\{s(t)=s \}}
- \sum_{s\in\mathcal{S}} \ind_{\{s(t)=s \}}\frac{1}{2}.\frac{\sum_{s' \in \mathcal{S}}V_i(t,s'){b}_{1iss'} }{r_i+ \sum_{s' \in \mathcal{S}}  V_i(t,s')b_{2iss'}} \\
-  \sum_{s\in\mathcal{S}} \ind_{\{s(t)=s \}}\frac{1}{2}.\frac{\epsilon_i+\sum_{s' \in \mathcal{S}}V_i(t,s')\bar{b}_{1iss'} }{\bar{r}_i+ \sum_{s' \in \mathcal{S}}  V_i(t,s')\bar{b}_{2iss'}} ,
\end{array}
 \end{equation}
and the equilibrium cost is $V_i(t,s),$ which satisfies the following ordinary differential system:

\begin{equation}
\begin{array}{ll}
\dot{V}_i  + \sum_{s' \in \mathcal{S}}V_i(t,s')\bar{b}_{ojss'} \\
-\frac{1}{4}\frac{(\sum_{s' \in \mathcal{S}}V_i(t,s'){b}_{1iss'})^2 }{[r_i+ \sum_{s' \in \mathcal{S}}  V_i(t,s')b_{2iss'}]}
-\frac{1}{4}\frac{(\epsilon_i+\sum_{s' \in \mathcal{S}}V_i(t,s')\bar{b}_{1iss'})^2 }{[\bar{r}_i+ \sum_{s' \in \mathcal{S}}  V_i(t,s')\bar{b}_{2iss'}]}\\
+\frac{1}{4}\sum_{j\neq i}[ \sum_{s' \in \mathcal{S}}  V_i(t,s')b_{2jss'} ]  (\frac{\sum_{s' \in \mathcal{S}}V_j(t,s'){b}_{1jss'} }{r_j+ \sum_{s' \in \mathcal{S}}  V_j(t,s')b_{2jss'}})^2\\
+\frac{1}{4}\sum_{j\neq i}[ \sum_{s' \in \mathcal{S}}  V_i(t,s')\bar{b}_{2jss'} ]  (\frac{\epsilon_j+\sum_{s' \in \mathcal{S}}V_j(t,s')\bar{b}_{1jss'} }{\bar{r}_j+ \sum_{s' \in \mathcal{S}}  V_j(t,s')\bar{b}_{2jss'}})^2\\

-\frac{1}{2}\sum_{j\neq i}[ \sum_{s' \in \mathcal{S}}  V_i(t,s')b_{1jss'} ]  (\frac{\sum_{s' \in \mathcal{S}}V_j(t,s'){b}_{1jss'} }{r_j+ \sum_{s' \in \mathcal{S}}  V_j(t,s')b_{2jss'}})\\
-\frac{1}{2}\sum_{j\neq i}[ \sum_{s' \in \mathcal{S}}  V_i(t,s')\bar{b}_{1jss'} ]  (\frac{\epsilon_j+\sum_{s' \in \mathcal{S}}V_j(t,s')\bar{b}_{1jss'} }{\bar{r}_j+ \sum_{s' \in \mathcal{S}}  V_j(t,s')\bar{b}_{2jss'}})=0,\\

{V}_i(T,s)=q_{i}(T,s), \ s\in \mathcal{S},
\end{array}
 \end{equation}
 \eod
\end{prop}

\textbf{Proof.} 

Step 1: We observe that the structure of the problem does not have a drift and is driven by regime switching.

Step 2: Based on step 1, we propose guess in the following form: $V_i(t,s).$ 

Step 3: we use the  stochastic integration formula for regime-switching to compute the difference between the cost and the guess functional as:
\begin{equation}
\begin{array}{ll}
\mathbb{E}  [L_i(s,u)-V_i(0,s_0)]
%=
%\mathbb{E} [q_{iT}- V_i(T, s(T))]\\
%+\mathbb{E}  \int_{0}^{T} \dot{V}_i  + \sum_{s' \in \mathcal{S}}V_i(t,s')\bar{b}_{ojss'} dt\\
%+\mathbb{E}  \int_{0}^{T}r_i (u_i-\bar{u}_i)^2+\sum_{j}[ \sum_{s' \in \mathcal{S}}  V_i(t,s')b_{2jss'} ](u_j-\bar{u}_j)^2 
%+\sum_{j}[ \sum_{s' \in \mathcal{S}}V_i(t,s'){b}_{1jss'}](u_j -\bar{u}_j)  dt\\
%+\mathbb{E}  \int_{0}^{T}\bar{r}_i\bar{u}_i^2+\sum_{j}[ \sum_{s' \in \mathcal{S}}V_i(t,s')\bar{b}_{2jss'}]\bar{u}_j^2 
%+ \epsilon_{i}\bar{u}_i+\sum_{j} [ \sum_{s' \in \mathcal{S}}V_i(t,s')\bar{b}_{1jss'}]\bar{u}_j dt\\
%%%
=
\mathbb{E} [q_{iT}- V_i(T, s(T))]\\
+\mathbb{E}  \int_{0}^{T} \dot{V}_i  + \sum_{s' \in \mathcal{S}}V_i(t,s')\bar{b}_{ojss'} \\
-\frac{1}{4}\frac{(\sum_{s' \in \mathcal{S}}V_i(t,s'){b}_{1iss'})^2 }{r_i+ \sum_{s' \in \mathcal{S}}  V_i(t,s')b_{2iss'}}\\
+\frac{1}{4}\sum_{j\neq i}[ \sum_{s' \in \mathcal{S}}  V_i(t,s')b_{2jss'} ]  (\frac{\sum_{s' \in \mathcal{S}}V_j(t,s'){b}_{1jss'} }{r_j+ \sum_{s' \in \mathcal{S}}  V_j(t,s')b_{2jss'}})^2\\
-\frac{1}{2}\sum_{j\neq i}[ \sum_{s' \in \mathcal{S}}  V_i(t,s')b_{1jss'} ]  (\frac{\sum_{s' \in \mathcal{S}}V_j(t,s'){b}_{1jss'} }{r_j+ \sum_{s' \in \mathcal{S}}  V_j(t,s')b_{2jss'}})\\
-\frac{1}{4}\frac{(\epsilon_i+\sum_{s' \in \mathcal{S}}V_i(t,s')\bar{b}_{1iss'})^2 }{[\bar{r}_i+ \sum_{s' \in \mathcal{S}}  V_i(t,s')\bar{b}_{2iss'}]},\\
+\frac{1}{4}\sum_{j\neq i}[ \sum_{s' \in \mathcal{S}}  V_i(t,s')\bar{b}_{2jss'} ]  (\frac{\epsilon_j+\sum_{s' \in \mathcal{S}}V_j(t,s')\bar{b}_{1jss'} }{\bar{r}_j+ \sum_{s' \in \mathcal{S}}  V_j(t,s')\bar{b}_{2jss'}})^2\\
-\frac{1}{2}\sum_{j\neq i}[ \sum_{s' \in \mathcal{S}}  V_i(t,s')\bar{b}_{1jss'} ]  (\frac{\epsilon_j+\sum_{s' \in \mathcal{S}}V_j(t,s')\bar{b}_{1jss'} }{\bar{r}_j+ \sum_{s' \in \mathcal{S}}  V_j(t,s')\bar{b}_{2jss'}})\\
+[r_i+ \sum_{s' \in \mathcal{S}}  V_i(t,s')b_{2iss'}] \{  u_i-\bar{u}_i+ \frac{1}{2}\frac{\sum_{s' \in \mathcal{S}}V_i(t,s'){b}_{1iss'} }{r_i+ \sum_{s' \in \mathcal{S}}  V_i(t,s')b_{2iss'}}  \}^2 dt\\
 +[\bar{r}_i+ \sum_{s' \in \mathcal{S}}  V_i(t,s')\bar{b}_{2iss'}] \{  \bar{u}_i+ \frac{1}{2}\frac{\epsilon_i+\sum_{s' \in \mathcal{S}}V_i(t,s')\bar{b}_{1iss'} }{\bar{r}_i+ \sum_{s' \in \mathcal{S}}  V_i(t,s')\bar{b}_{2iss'}}  \}^2dt
\end{array}
 \end{equation}

Step 4: Assuming  that $ [r_i+ \sum_{s' \in \mathcal{S}}  V_i(t,s')b_{2iss'}] >0,$ and $[\bar{r}_i+ \sum_{s' \in \mathcal{S}}  V_i(t,s')\bar{b}_{2iss'}] >0$ 
the terms completion lead to a one-shot optimization of a strictly convex and coercive function.

Step 5: the minimization  and the identification of the processes provides the announced result.

 }  

\subsection{Quadratic-Quadratic  MFTG}
This example examines a class of Quadratic-Quadratic Mean-Field-Type Game (QQ-MFTG)  problem. The state is non-linear in $u.$ A semi-explicit solution is derived.
\begin{equation} \label{eq:QQ}
\begin{cases}
L_i(x,u) = q_{iT} x_T+\int_0^T r_i (u_i-\bar{u}_i)^2+\bar{r}_i\bar{u}_i^2 + \bar{\epsilon}_{2i}\bar{u}_i dt\\
\inf_{u_i} ~\mathbb{E}[L_i(x,u)],\\ 
\mbox{subject to }\\
dx= \sum_{j \in \mathcal{I}} [q_j (u_j-\bar{u}_j)^2+\bar{q}_j\bar{u}_j^2 + \epsilon_{1j}(u_j - \bar{u}_j) + \bar{\epsilon}_{1j}\bar{u}_j]dt\\ +\sigma dB+\int \mu  d\tilde{N},\\
x(0)=x_0.
\end{cases}
\end{equation} 
where the coefficients are regime-switching dependent with switching rate matrix $\tilde{Q}=(\tilde{q}_{ss'}, (s,s')\in \mathcal{S}^2).$ 

\begin{prop} \label{semi:prop5}
 Assume that $r_i>0, q_i\geq 0,\ r_i>\delta>0, \bar{r}_i>0, \bar{q}_i\geq 0, \bar{r}_i>\delta >0,$ $ \mathbb{E}[\epsilon_{1i}] = 0,~\forall i \in \mathcal{I}.$
The  QQ-MFTG problem   (\ref{eq:QQ})  has unique solution and it is given by
\begin{equation}
\begin{array}{ll}
u_i^*-\bar{u}_i^* =  -\sum_{s\in\mathcal{S}} \ind_{\{s(t)=s \}}  \frac{\alpha_i}{2(r_i+ \alpha_i  q_i)}  \epsilon_{1i},\\
\bar{u}_i^* = - \sum_{s\in\mathcal{S}} \ind_{\{s(t)=s \}} [\frac{\left[ \bar{\epsilon}_{2i} + \alpha_i \bar{\epsilon}_{1i}\right]}{2 \left[\bar{r}_i + \alpha_i \bar{q}_i \right]}],\\
\mathbb{E}L_i^* = \mathbb{E}[\alpha_i(0,s(0))x_0] \\
+ \mathbb{E} \int_{0}^{T}\left\{ \right. \\
 \alpha_i \sum_{j \ne i} q_j \left(\frac{\alpha_j}{2(r_j+ \alpha_j  q_j)}  \epsilon_{1j} \right)^2+ \alpha_i \sum_{j \ne i} \bar{q}_j \left(\frac{\left[ \bar{\epsilon}_{2j} + \alpha_j \bar{\epsilon}_{1j}\right]}{2 \left[\bar{r}_j + \alpha_j \bar{q}_j \right]}\right)^2   \\
 -  \left( \alpha_i \sum_{j \ne i} \epsilon_{1j}\left(\frac{\alpha_j}{2(r_j+ \alpha_j  q_j)}  \epsilon_{1j}\right) + \alpha_i \sum_{j \ne i} \bar{\epsilon}_{1j} \left( \frac{\left[ \bar{\epsilon}_{2j} + \alpha_j \bar{\epsilon}_{1j}\right]}{2 \left[\bar{r}_j + \alpha_j \bar{q}_j \right]} \right)  \right)  \\
 - \frac{\alpha_i^2}{4(r_i+ \alpha_i  q_i)}  \epsilon_{1i}^2 
 -  \frac{\left[ \bar{\epsilon}_{2i} + \alpha_i \bar{\epsilon}_{1i}\right]^2}{4 \left[\bar{r}_i + \alpha_i \bar{q}_i \right]}  \left. \right\} dt\\
\end{array}
\end{equation} 
Note that the semi-explicit solution is in fact an explicit solution. Let $\vec{\alpha}_i(t) = [\alpha_i(t,s)]_{s \in \mathcal{S}}$, and $\vec{q}_i(T) = [q_i(T,s)]_{s \in \mathcal{S}}$, then $\vec{\alpha}_i(t)$ is explicitly given by
\begin{align*}
\vec{\alpha}_i(t) = \vec{q}_i(T) \exp\left[{-\int_{t}^{T} \tilde{Q} \dif t'}\right].
\end{align*}
in particular $\vec{\alpha}_i(0) = \vec{q}_i(T) \exp[{-\int_{0}^{T} \tilde{Q} \dif t'}]$. \eod\\
\end{prop}

\textbf{Proof.} Let us consider the following guess functional: $V_i = \alpha_i x$. Then,
\begin{align*}
&\partial_t V_i = \dot{\alpha}_i x,\ 
\partial_x V_i = \alpha_i,\\
&V_i(t,x+\mu) - V_i(t,x) - \partial_x V_i \mu  = 0,
\end{align*}
It\^o's formula yields
\begin{equation} \begin{array}{ll}
V_i(T) - V_i(0) = \int_{0}^{T} \left( \dot{\alpha}_i x + \alpha_i \sum_{j \in \mathcal{I}} q_j (u_j-\bar{u}_j)^2 \right. \\ 
+ \alpha_i \sum_{j \in \mathcal{I}} \bar{q}_j\bar{u}_j^2 + \alpha_i \sum_{j \in \mathcal{I}} \epsilon_{1j}(u_j - \bar{u}_j)\\
 + \alpha_i \sum_{j \in \mathcal{I}} \bar{\epsilon}_{1j}\bar{u}_j 
+ \sum_{s' \in \mathcal{S}} x \left( \alpha_i(.,s') - \alpha_i(.,s)  \right) \tilde{q}_{ss'} \left.  \right)  dt \\
+ \int_{0}^{T}  \alpha_i \sigma \dB + \int_{\Theta} \alpha_i \mu \tilde{N} (dt,\dif \theta).
\end{array}
\end{equation}
Thus, the difference $\mathbb{E}[L_i - V_i(0)]$ is given by
\begin{equation} 
\begin{array}{ll}
\mathbb{E}[L_i - V_i(0)] = \mathbb{E} \left[ q_{i}(T) x(T) - V_i(T) \right] \\
+ \mathbb{E} \int_0^T  r_i (u_i-\bar{u}_i)^2+\bar{r}_i\bar{u}_i^2 + \bar{\epsilon}_{2i}\bar{u}_i \\
+
 \dot{\alpha}_i x + \alpha_i \sum_{j \in \mathcal{I}} q_j (u_j-\bar{u}_j)^2  \\
+ \alpha_i \sum_{j \in \mathcal{I}} \bar{q}_j\bar{u}_j^2 + \alpha_i \sum_{j \in \mathcal{I}} \epsilon_{1j}(u_j - \bar{u}_j)\\
 + \alpha_i \sum_{j \in \mathcal{I}} \bar{\epsilon}_{1j}\bar{u}_j   
+ x\sum_{s' \in \mathcal{S}}  \left( \alpha_i(.,s') - \alpha_i(.,s)  \right) \tilde{q}_{ss'} \ dt
\end{array}
\end{equation}
Performing square completion  yields  
\begin{equation} \begin{array}{ll}
\left( (u_i-\bar{u}_i)^2 + \frac{\alpha_i}{r_i+ \alpha_i  q_i}  \epsilon_{1i}(u_i - \bar{u}_i) \right) \\
= \left(u_i-\bar{u}_i + \frac{\alpha_i}{2(r_i+ \alpha_i  q_i)}  \epsilon_{1i} \right)^2 - \frac{\alpha_i^2}{4(r_i+ \alpha_i  q_i)^2}  \epsilon_{1i}^2,\\
\left(  \bar{u}_i^2 + \frac{\left[ \bar{\epsilon}_{2i} + \alpha_i \bar{\epsilon}_{1i}\right]}{\left[\bar{r}_i + \alpha_i \bar{q}_i \right]} \bar{u}_i  \right) \\
= \left( \bar{u}_i + \frac{\left[ \bar{\epsilon}_{2i} + \alpha_i \bar{\epsilon}_{1i}\right]}{2 \left[\bar{r}_i + \alpha_i \bar{q}_i \right]}  \right)^2
 - \frac{\left[ \bar{\epsilon}_{2i} + \alpha_i \bar{\epsilon}_{1i}\right]^2}{4 \left[\bar{r}_i + \alpha_i \bar{q}_i \right]^2},
\end{array}
\end{equation}

\begin{equation}
\begin{array}{ll}
\mathbb{E}[L_i - V_i(0)] =\mathbb{E} \left[ q_{i}(T) x(T) - V_i(T) \right]  \\
+ \mathbb{E} \int_{0}^{T}
 \dot{\alpha}_i x + \alpha_i \sum_{j \ne i} q_j (u_j-\bar{u}_j)^2+ \alpha_i \sum_{j \ne i} \bar{q}_j\bar{u}_j^2 \\
+ \alpha_i \sum_{j \ne i} \epsilon_{1j}(u_j - \bar{u}_j) + \alpha_i \sum_{j \ne i} \bar{\epsilon}_{1j}\bar{u}_j   \\
+  (r_i+ \alpha_i  q_i) \left(u_i-\bar{u}_i + \frac{\alpha_i}{2(r_i+ \alpha_i  q_i)}  \epsilon_{1i} \right)^2 \\
 -  \frac{\alpha_i^2}{4(r_i+ \alpha_i  q_i)}  \epsilon_{1i}^2 \\
 +  \left[\bar{r}_i + \alpha_i \bar{q}_i \right] \left( \bar{u}_i + \frac{\left[ \bar{\epsilon}_{2i} + \alpha_i \bar{\epsilon}_{1i}\right]}{2 \left[\bar{r}_i + \alpha_i \bar{q}_i \right]}  \right)^2  \\
 - \frac{\left[ \bar{\epsilon}_{2i} + \alpha_i \bar{\epsilon}_{1i}\right]^2}{4 \left[\bar{r}_i + \alpha_i \bar{q}_i \right]} 
+ \sum_{s' \in \mathcal{S}} x \left( \alpha_i(.,s') - \alpha_i(.,s)  \right) \tilde{q}_{ss'}  
\end{array}
\end{equation}

Minimizing terms it yields
\begin{align*}
\dot{\alpha}_i(t,s) &= -  \sum_{s' \in \mathcal{S}} \left( \alpha_i(.,s') - \alpha_i(.,s)  \right) \tilde{q}_{ss'} ,\\
\alpha_i(T,s) &= q_i(T,s),
\end{align*}
completing the proof. \eod 

\begin{remark}
 Notice that using the result presented in Proposition \ref{semi:prop5}, the following Quadratic-Exponential-Quadratic Mean-Field-Type Game (QEQ-MFTG) problem: 
 \begin{equation}
 \begin{cases}
L_i(x,u) = q_{iT} x_T+\int_0^T r_i (u_i-\bar{u}_i)^2+\bar{r}_i\bar{u}_i^2 + \bar{\epsilon}_{2i}\bar{u}_i  dt,\\
 \inf_{u_i}  \mathbb{E} \exp \left[{\lambda_i L_i(x,u)} \right],\\
 \mbox{subject to }\\
 dx= \sum_{j \in \mathcal{I}} [q_j (u_j-\bar{u}_j)^2+\bar{q}_j\bar{u}_j^2 + \epsilon_{1j}(u_j - \bar{u}_j) + \bar{\epsilon}_{1j}\bar{u}_j]dt\\ +\sigma dB+\int \mu  d\tilde{N},\\
 x(0)=x_0.
 \end{cases}
 \end{equation}
 can be solved explicitly. \eod
\end{remark}

\subsection{Quadratic State and Power Utility }
{\color{black} This subsection examines a class of mean-field-type games with power payoffs and a non-linear state. 
 The model is inspired from the modern portfolio optimization under shared asset platform by several decision-makers. The state $x(t)$ is the total amount of money. Decision-maker $i$ can decide to consume certain amount $u_{1i}$ and re-allocate the remaining between less-risky assets $(1-u_{2i})\kappa_2 x$ and more risky assets $u_{2i}\kappa_2 x + 
 \bar{u}_{2i}{x} [{\sigma} dB+\int {\mu}  d\tilde{N}].$  The  coefficients $\kappa_1, \kappa_2$  depend on time $t$ and on the switching regime $s(t)$ which takes values in $\mathcal{S}.$  The set $S$ is non-empty and finite.
 We have modified the model to include mean-field terms, a function of the expected value of the state and a function of the expected value of the control action.  }
\begin{equation} \label{coobdouglas1}
\begin{cases}
R_i(x,u) =-q_{iT} \dfrac{(x_T-\bar{x}_T)^{2k_i}}{2k_i}+ \bar{q}_{iT} \dfrac{\bar{x}^{\bar{\rho}_i}_T}{\bar{\rho}_i}\\ 
+\int_0^T -q_i \dfrac{(x-\bar{x})^{2k_i}}{2k_i} +\bar{q}_i\dfrac{\bar{x}^{\bar{\rho}_i}}{\bar{\rho}_i} 
- r_i \dfrac{(u_{1i}-\bar{u}_{1i})^{2k_i}}{2k_i} +\bar{r}_i\dfrac{\bar{u}_{1i}^{\bar{\rho}_i}}{\bar{\rho}_i} dt,\\
\sup_{u_i}  \mathbb{E}[R_i(x,u)],\\
\mbox{subject to }\\
dx= \sum_{i \in \mathcal{I}}[-(u_{1i}-\bar{u}_{1i})+(1- (u_{2i}-\bar{u}_{2i}))\kappa_1 (x-\bar{x})]\\
+\sum_{i \in \mathcal{I}}  (u_{2i}-\bar{u}_{2i})\kappa_2 (x-\bar{x})dt\\
 +\sum_{i \in \mathcal{I}} (u_{2i}-\bar{u}_{2i}) (x-\bar{x}) [\sigma dB+\int \mu  d\tilde{N}]\\
+ \sum_{i \in \mathcal{I}} [-\bar{u}_{1i}+(1- \bar{u}_{2i})\bar{\kappa}_1 \bar{x}+ \bar{u}_{2i}\bar{\kappa}_2 \bar{x}]dt\\
 +\sum_{i \in \mathcal{I}} \bar{u}_{2i}\bar{x} [\bar{\sigma} dB_o+\int \bar{\mu}  d\tilde{N}_o],\\
x(0)=x_0.
\end{cases} 
\end{equation}
where the coefficients are regime-switching dependent, and $k_i\geq 1,  \bar{\rho}_i\in (0,1).$ The coefficients $q_i,\bar{q}_i, r_i,\bar{r}_i$ are positive. The state dynamics (\ref{coobdouglas1}) is not linear in $(x,u).$  %%because the term $x u$ and $\bar{x}\bar{u}$ are non-linear.

 {Following the same method as in the problem \eqref{coobdouglas1}, a semi-explicit solution can be derived.}

Note that a similar method can be used to derive semi-explicit solution to the following game problem in which decision-makers minimize with $k_i(t,s)> 1, \  \bar{\rho}_i(t,s)>1.$

\begin{equation}
\begin{cases}
L_i(x,u) = q_{iT} \frac{(x_T-\bar{x}_T)^{2 k_i}}{2k_i}+ \bar{q}_{iT} \frac{\bar{x}^{\bar{\rho_i}}_T}{\bar{\rho}_i} \\ +\int_0^T q_i \frac{(x-\bar{x})^{2k_i}}{2k_i} +\bar{q}_i\frac{\bar{x}^{\bar{\rho}_i}}{\bar{\rho}_i} 
+ r_i \frac{(u_{1i}-\bar{u}_{1i})^{2k_i}}{2k_i} +\bar{r}_i\frac{\bar{u}_{1i}^{\bar{\rho}_i}}{\bar{\rho}_i} dt,\\
\inf_{u_i} \mathbb{E}[L_i(x,u)],\\
\mbox{subject to } \\
dx=\sum_{i} (u_{1i}-\bar{u}_{1i})+(1- (u_{2i}-\bar{u}_{2i}))k_1 (x-\bar{x})\\
+ \sum_{i}(u_{2i}-\bar{u}_{2i})k_2 (x-\bar{x})dt\\
 +\sum_{i}(u_{2i}-\bar{u}_{2i}) (x-\bar{x}) [\sigma dB+\int \mu  d\tilde{N}]\\
+\sum_{i} [-\bar{u}_{1i}+(1- \bar{u}_{2i})\bar{k}_1 \bar{x}+ \bar{u}_{2i}\bar{k}_2 \bar{x}]dt\\
 +\sum_{i}\bar{u}_{2i}\bar{x} [\bar{\sigma} dB_o+\int \bar{\mu}  d\tilde{N}_o],\\
x(0)=x_0.
\end{cases}
\end{equation}

This can be easily extended to include multi-type power utilities in the following form:

\begin{equation}
\begin{array}{ll}
\sum_{k=1}^K q_{ik} \frac{(x-\bar{x})^{\rho_{ik}}}{\rho_{ik}} +\bar{q}_{ik}\frac{\bar{x}^{\bar{\rho}_{ik}}}{\bar{\rho}_{ik}} \\ 
+\sum_{k=1}^K r_{ik} \frac{(u_{1i}-\bar{u}_{1i})^{\rho_{ik}}}{\rho_{ik}} +\bar{r}_{ik}\frac{\bar{u}^{\bar{\rho}_{ik}}_{1i}}{\bar{\rho}_{ik}}. 
\end{array}
\end{equation}
In this case, the  guess functional will be
$$
\sum_{k=1}^K p_{ik} \frac{(x-\bar{x})^{\rho_{ik}}}{\rho_{ik}} +\bar{p}_{ik}\frac{\bar{x}^{\bar{\rho}_{ik}}}{\bar{\rho}_{ik}},
$$
with an ordinary differential system for $p_{ik},$ and $ \bar{p}_{ik}.$

\subsection{Non-Linear State and Log-Utility }
We consider the following logarithmic   Cobb-Douglas  utility.
\begin{equation} \label{coobdouglas2}
\begin{cases}
R_i(x,u) =-q_{iT} \frac{(x_T-\bar{x}_T)^{2k_i}}{2k_i} +\log [ ( \bar{x}_T)^{\bar{q}_{iT} }] \\ +\int_0^T 
-q_{i} \frac{(x-\bar{x})^{2k_i}}{2k_i}- r_i\frac{(u_{1i}-\bar{u}_{1i})^{2k_i}}{2k_i}+
  \log[ \bar{x}^{\bar{q}_i}\bar{u}_{1i}^{\bar{r}_i}] dt\\
\sup_{u_i} \mathbb{E} [ R_i(x,u) ],\\
\mbox{subject to }\\
dx= \sum_{i}[-(u_{1i}-\bar{u}_{1i})+(1- (u_{2i}-\bar{u}_{2i}))\kappa_1 (x-\bar{x})]\\
+\sum_{i} (u_{2i}-\bar{u}_{2i})\kappa_2 (x-\bar{x})dt\\
 +\sum_{i}(u_{2i}-\bar{u}_{2i}) (x-\bar{x}) [\sigma dB+\int \mu  d\tilde{N}]\\
+\sum_{i} [-\bar{u}_{1i}+(1- \bar{u}_{2i})\bar{\kappa}_1 \bar{x}+ \bar{u}_{2i}\bar{\kappa}_2 \bar{x}]dt\\
 +\sum_{i}\bar{u}_{2i}\bar{x} [\bar{\sigma} dB_o+\int \bar{\mu}  d\tilde{N}_o],\\
x(0)=x_0,\ s(0)=s_0,\\
\mathbb{P}(s(t+\epsilon)=s' | s,u)=\int_t^{t+\epsilon} \tilde{q}_{ss'} dt' + o(\epsilon),\ s'\neq s
\end{cases}
\end{equation}
where the coefficients are regime-switching dependent.

Note that the state dynamics (\ref{coobdouglas2}) is not linear in $(x,u).$

{Following the same method as above, the problem \eqref{coobdouglas2} can be solved explicitly.}

%\subsection{Non-Linear State and Exponential-Utility }
%
%\begin{equation}
%\begin{cases}
%L_i(x,u) = q_{iT} (x_T-\bar{x}_T)^2+ \bar{q}_{iT}e^{\bar{x}_T}\\
%+\int_0^T q_{i} ({x-\bar{x}})^2+ \bar{q}_{i}e^{\bar{x}}+
%r_i (u_i-\bar{u}_i)^2+ e^{\bar{x}}[\bar{r}_i\bar{u}_i^2 + \epsilon_{1i}\bar{u}_i] dt,\\
%\inf_{u_i} \mathbb{E}[ L_i(x,u) ],\\
%\mbox{subject to }\\
%dx= \sum_{j}[ b_{2j} (u_j-\bar{u}_j)+  \bar{b}_{2j} \bar{u}_j]dt
%+\sigma dB+\int \mu  d\tilde{N},\\
%x(0)=x_0.
%\end{cases}
%\end{equation}
%where the coefficients are regime-switching dependent.
%
% \textcolor{blue}{Following the same method as above, the problem \eqref{coobdouglas2} can be solved explicitly.}

\subsection{Cotangent  Drift}
This subsection we examine mean-field-type games with cotangent  drift. This class of games is inspired from \cite{refpp1,refpp2,refpp3, refpp4, refpp5,refpp6,refpp7,refpp8,refpp9,refpp10}. We have modified the model to include mean-field terms. Using trigonometric relationships a semi-explicit equilibrium solution is derived. 
\begin{equation} \label{eq:cotangent}
\begin{array}{ll} 
L_i(x,u) = \\
\int_0^T ((u_i-\bar{u}_i)^2-q_i)\cos^2(\frac{x-\bar{x}}{4})+q_i  +
(\bar{u}_i^2-\bar{q}_i)\cos^2(\frac{\bar{x}}{4})+\bar{q}_i dt,\\
\inf_{u_i} \mathbb{E}[ L_i(x,u) ],\\
\mbox{subject to }\\
dx= [\frac{1}{2}cot(\frac{x-\bar{x}}{2})+\sum_{j}b_{2j}(u_j-\bar{u}_j)+
\frac{1}{2}cot(\frac{\bar{x}}{2})+\sum_{j}\bar{b}_{2j} \bar{u}_j]dt\\
+\sigma dB,\\
x(0)=x_0.
\end{array}
\end{equation}
where the coefficients are regime-switching dependent and $cot(\theta)=\frac{1}{\tan(\theta)}=-cot(-\theta).$ 
\begin{prop} \label{semi:prop6}
Assume that $q_i>0, \bar{q}_i>0.$
The mean-field-type game  problem with cotangent drift (\ref{eq:cotangent}) has a unique equilibrium  solution which is given by
\begin{equation} \label{eq:cotangentsol1}
\begin{array}{ll} 
u^*_i-\bar{u}^*_i= - \sum_{s\in\mathcal{S}} \ind_{\{s(t)=s \}} b_{2i}\frac{\alpha_i}{4} tan(\frac{x-\bar{x}}{4})\\
\bar{u}^*= -\sum_{s\in\mathcal{S}} \ind_{\{s(t)=s \}} \frac{\bar{\alpha}_i}{4} \bar{b}_{2i}tan(\frac{\bar{x}}{4})\\
\mathbb{E}[ L_i(x,u^*) ]= \\ \mathbb{E}[\alpha_i(0,s_0) \sin^2(\frac{x_0-\bar{x}_0}{4}) + \bar{\alpha}_i(0,s_0)  \sin^2(\frac{\bar{x}_0}{4})+\delta_i(0,s_0) ],
 \\
\end{array}
\end{equation}
whenever the following system

\begin{equation} \label{eq:cotangentsol2}
\begin{array}{ll} 
\dot{\alpha}_i +  q_i -(2+\sigma^2)\frac{\alpha_i}{8} - \frac{ \alpha_i^2}{16} b_{2i}^2 \\
+ \sum_{s'}[\alpha_i(t,s')-\alpha_i(t,s)]\tilde{q}_{ss'}
-\frac{\alpha_i}{8} \sum_{j\neq i}\alpha_j b_j^2 =0,\\

{\alpha}_i(T,s)=0,\\
%\dot{\bar{\alpha}}_i \\
 \dot{\bar{\alpha}}_i+\bar{q}_i - \frac{\bar{\alpha}_i}{4}- \frac{\bar{\alpha}^2_i}{16} \bar{b}^2_{2i}\\ 
+ \sum_{s'}[\bar{\alpha}_i(t,s')-\bar{\alpha}_i(t,s)]\tilde{q}_{ss'} - \frac{\bar{\alpha}_i}{8} \sum_{j\neq i}\bar{b}^2_{2j} \bar{\alpha}_j =0,\\
\bar{\alpha}_i(T,s)=0,\\
%\dot{\delta}_i+ ...  \\
\dot{\delta}_i+\frac{\alpha_i}{16} (2+\sigma^2)+ \frac{\bar{\alpha}_i}{8}+
\sum_{s'}[\delta_i(t,s')-\delta_i(t,s)]\tilde{q}_{ss'} =0,\\
\delta_i(T,s)=0,
\end{array} 
\end{equation} 
has a unique solution  with positive $\alpha_i, \bar{\alpha}_i$ which do not blow up within $[0,T].$ \eod
\end{prop}

{\bf Proof}: \\
We prove the statement using a direct method.  Step 1: We observe that
the mean-field-type problem is driven by functionals of $x-\bar{x}$ and $\bar{x}$ which are conditionally orthogonal processes.

 Step 2:  Given the structure of the problem, we propose the following guess functional:
 $$f_i= {\alpha}_i(t,s) \sin^2(\frac{x-\bar{x}}{4}) + \bar{\alpha}_i(t,s)  \sin^2(\frac{\bar{x}}{4})+\delta_i(t,s),$$ be a guess functional. 

Step 3: we apply Brownian with regime switching to obtain the difference between the cost functional and the guess functional as:
\begin{equation} \label{eq:cotangentsol3}
\begin{array}{ll} 
\mathbb{E}[ L_i -f_i(0)]=  -\mathbb{E}f_i(T)\\
+\mathbb{E}\int_0^T  \cos^2(\frac{x-\bar{x}}{4}) \left(u-\bar{u} + b_{2i}\frac{\alpha_i}{4} tan(\frac{x-\bar{x}}{4})\right)^2 dt\\
+\mathbb{E}\int_0^T  \sin^2(\frac{x-\bar{x}}{4}) \left\{  \dot{\alpha}_i +  q_i -(2+\sigma^2)\frac{\alpha_i}{8} - \frac{ \alpha_i^2}{16} b_{2i}^2 \right. \\ \left. 
+ \sum_{s'}[\alpha_i(t,s')-\alpha_i(t,s)]\tilde{q}_{ss'}
-\frac{\alpha_i}{8} \sum_{j\neq i}\alpha_j b_j^2 \right\} dt\\
+\dot{\delta}_i+\frac{\alpha_i}{16} (2+\sigma^2)+
\sum_{s'}[\delta_i(t,s')-\delta_i(t,s)]\tilde{q}_{ss'}\\
%% \dot{\bar{\alpha}}_i  \sin^2(\frac{\bar{x}}{4})+ \frac{\bar{\alpha}_i}{4} \sin(\frac{\bar{x}}{2})
%%[\frac{1}{2}cot(\frac{\bar{x}}{2})+\sum_{j}\bar{b}_{2j} \bar{u}_j]]
%%%+(\bar{u}_i^2-\bar{q}_i)\cos^2(\frac{\bar{x}}{4})+\bar{q}_i dt,

+\left\{ \dot{\bar{\alpha}}_i+\bar{q}_i - \frac{\bar{\alpha}_i}{4}- \frac{\bar{\alpha}^2_i}{16} \bar{b}^2_{2i} \right. \\ 
+ \sum_{s'}[\bar{\alpha}_i(t,s')-\bar{\alpha}_i(t,s)]\tilde{q}_{ss'} - \frac{\bar{\alpha}_i}{8} \sum_{j\neq i}\bar{b}^2_{2j} \bar{\alpha}_j \left. \right\}  \sin^2(\frac{\bar{x}}{4}) dt\\
+
\cos^2(\frac{\bar{x}}{4})[  \bar{u}_i +  \frac{\bar{\alpha}_i}{4} \bar{b}_{2i}tan(\frac{\bar{x}}{4}) ]^2 
+\frac{\bar{\alpha}_i}{8}  dt,
\end{array}
\end{equation}

Step 4: Noting the terms completion leads to a strictly concave one-shot optimization with coercive  function
$\cos^2(\frac{\bar{x}}{4})[  \bar{u}_i +  \frac{\bar{\alpha}_i}{4} \bar{b}_{2i}tan(\frac{\bar{x}}{4}) ]^2$  whenever  $\cos^2(\frac{\bar{x}}{4})>0.$ 

Step 5:
By identification of processes one obtains the announced result. This completes the proof.
 \eod

{\color{blue}
The mean-field term $\bar{x}$ solves
\begin{equation} \label{eq:cotangentsol3t1yy}
\begin{array}{ll}
d\bar{x}= [\frac{1}{2}cot(\frac{\bar{x}}{2})-\sum_{j}\bar{b}^2_{2j} \frac{\bar{\alpha}_j}{4}\tan(\frac{\bar{x}}{4})]dt,\\
\bar{x}(0)=\bar{x}_0,
\end{array}
\end{equation}
which has a unique solution for $\bar{x}_0\in (0, \pi).$
}

\subsection{Hyperbolic coTangent  Drift}
Problem (\ref{eq:cotangent}) can be modified to handle the  hyperbolic cotangent  drift case as specified below. The functions $\cos, \sin, \tan, \cot$ are replaced by $\cosh,\sinh, \tanh, \coth$ respectively.
\begin{equation} \label{hyper}
\begin{array}{ll}
L_i(x,u) =  
\int_0^T ((u_i-\bar{u}_i)^2+q_i)\cosh^2(\frac{x-\bar{x}}{4})-q_i \\ +
(\bar{u}_i^2+\bar{q}_i)\cosh^2(\frac{\bar{x}}{4})-\bar{q}_i dt,\\
\inf_{u_i} \mathbb{E}[ L_i(x,u) ],\\
\mbox{subject to }\\
dx= [\frac{1}{2}coth(\frac{x-\bar{x}}{2})+\sum_{j}b_{2j}(u_j-\bar{u}_j)\\ +
\frac{1}{2}coth(\frac{\bar{x}}{2})+\sum_{j}\bar{b}_{2j} \bar{u}_j]dt
+\sigma dB,\\
x(0)=x_0.
\end{array}
\end{equation}
where the coefficients are regime-switching dependent and $\coth(\theta)=\frac{e^{\theta}+e^{\theta}}{e^{\theta}-e^{-\theta}}=-\coth(-\theta).$
\begin{prop} \label{semi:prop7}
The equilibrium strategies and the equilibrium costs are given by 
\begin{equation} \label{eq:cotangentsol1ui}
\begin{array}{ll} 
u^*_i-\bar{u}^*_i= - \sum_{s\in\mathcal{S}} \ind_{\{s(t)=s \}} b_{2i}\frac{\alpha_i}{4} \tanh(\frac{x-\bar{x}}{4})\\
\bar{u}^*_i=- \sum_{s\in\mathcal{S}} \ind_{\{s(t)=s \}} \frac{\bar{\alpha}_i}{4} \bar{b}_{2i}\tanh(\frac{\bar{x}}{4}) \\
\mathbb{E}[ L_i(x,u^*) ]= \mathbb{E}[\alpha_i(0,s_0) \sinh^2(\frac{x_0-\bar{x}_0}{4})]\\
 +\mathbb{E}[ \bar{\alpha}_i(0,s_0)  \sinh^2(\frac{\bar{x}_0}{4})+\delta_i(0,s_0) ],
 \\
\end{array}
\end{equation}
whenever the following system: 

\begin{equation} \label{eq:cotangentsol267}
\begin{array}{ll} 
\dot{\alpha}_i +  q_i +(2+\sigma^2)\frac{\alpha_i}{8} - \frac{ \alpha_i^2}{16} b_{2i}^2 \\
+ \sum_{s'}[\alpha_i(t,s')-\alpha_i(t,s)]\tilde{q}_{ss'}
-\frac{\alpha_i}{8} \sum_{j\neq i}\alpha_j b_j^2 =0,\\. 

{\alpha}_i(T,s)=0,\\
%\dot{\bar{\alpha}}_i \\
 \dot{\bar{\alpha}}_i+\bar{q}_i +\frac{\bar{\alpha}_i}{4}- \frac{\bar{\alpha}^2_i}{16} \bar{b}^2_{2i}\\ 
+ \sum_{s'}[\bar{\alpha}_i(t,s')-\bar{\alpha}_i(t,s)]\tilde{q}_{ss'} - \frac{\bar{\alpha}_i}{8} \sum_{j\neq i}\bar{b}^2_{2j} \bar{\alpha}_j =0,\\
\bar{\alpha}_i(T,s)=0,\\
%\dot{\delta}_i+ ...  \\
\dot{\delta}_i+\frac{\alpha_i}{16} (2+\sigma^2)+ \frac{\bar{\alpha}_i}{8}+
\sum_{s'}[\delta_i(t,s')-\delta_i(t,s)]\tilde{q}_{ss'} =0\\
\delta_i(T,s)=0,
\end{array} 
\end{equation}
has unique solution  with positive $(\alpha_i, \bar{\alpha}_i)_i$ which do not blow up within $[0,T].$ \eod
\end{prop}

 The system in (\ref{eq:cotangentsol267})  shares some similarities with the system in (\ref{eq:cotangentsol2}) of Problem (\ref{eq:cotangent}). However, these two systems are different. In particular, the sign  of the terms $(2+\sigma^2)\frac{\alpha_i}{8},$ and $\frac{\bar{\alpha}_i}{4}$ have changed.

{\bf Proof.}

We prove the statement on the hyperbolic game using a direct method.  
Let 
\begin{equation}  \nonumber
\begin{array}{ll} 
f_i(t,x,s)= {\alpha}_i(t,s(t)) \sinh^2(\frac{x(t)-\bar{x}(t)}{4}) \\ + \bar{\alpha}_i(t,s(t))  \sinh^2(\frac{\bar{x}(t)}{4})+\delta_i(t,s(t)),
\end{array}
\end{equation}
 be a guess functional combining hyperbolic functions. 
\begin{equation} \label{eq:cotangentsol3}
\begin{array}{ll} 
\mathbb{E}[ L_i -f_i(0)]=  -\mathbb{E}f_i(T)\\
+\mathbb{E}\int_0^T  \cosh^2(\frac{x-\bar{x}}{4}) \left(u-\bar{u} + b_{2i}\frac{\alpha_i}{4} \tanh(\frac{x-\bar{x}}{4})\right)^2 dt\\
+\mathbb{E}\int_0^T  \sinh^2(\frac{x-\bar{x}}{4}) \left\{  \dot{\alpha}_i +  q_i +(2+\sigma^2)\frac{\alpha_i}{8} - \frac{ \alpha_i^2}{16} b_{2i}^2 \right. \\ \left. 
+ \sum_{s'}[\alpha_i(t,s')-\alpha_i(t,s)]\tilde{q}_{ss'}
-\frac{\alpha_i}{8} \sum_{j\neq i}\alpha_j b_j^2 \right\} dt\\
+\dot{\delta}_i+\frac{\alpha_i}{16} (2+\sigma^2)+ \frac{\bar{\alpha}_i}{8} +
\sum_{s'}[\delta_i(t,s')-\delta_i(t,s)]\tilde{q}_{ss'}\\
+\left\{ \dot{\bar{\alpha}}_i+\bar{q}_i + \frac{\bar{\alpha}_i}{4}- \frac{\bar{\alpha}^2_i}{16} \bar{b}^2_{2i} \right. \\ 
+ \sum_{s'}[\bar{\alpha}_i(t,s')-\bar{\alpha}_i(t,s)]\tilde{q}_{ss'} - \frac{\bar{\alpha}_i}{8} \sum_{j\neq i}\bar{b}^2_{2j} \bar{\alpha}_j \left. \right\}  \sinh^2(\frac{\bar{x}}{4}) dt\\
+
\cosh^2(\frac{\bar{x}}{4})[  \bar{u}_i +  \frac{\bar{\alpha}_i}{4} \bar{b}_{2i}\tanh(\frac{\bar{x}}{4}) ]^2 
dt,
\end{array}
\end{equation}

By identification one obtains the announced result. This completes the proof.  \eod 

{\color{blue}
The mean-field term $\bar{x}$ solves
\begin{equation} \label{eq:cotangentsol3t1}
\begin{array}{ll}
d\bar{x}= [\frac{1}{2}coth(\frac{\bar{x}}{2})-\sum_{j}\bar{b}^2_{2j} \frac{\bar{\alpha}_j}{4}\tanh(\frac{\bar{x}}{4})]dt,\\
\bar{x}(0)=\bar{x}_0\neq 0,
\end{array}
\end{equation}
which has a unique global  solution within $[0,T].$
}

\subsection{ A Delayed and Trend-based MFTG}
We present a cooperative MFTG with basic state dynamics $x(t)$, regime switching $s(t),$ a trend $ y(t):=\int_{-\tau}^0 e^{\lambda t'}x(t+t') dt'$ on the time window $[t-\tau, t],$ the delayed state $z(t)=x(t-\tau).$ This class of examples plays  an important role in real-world applications as the effects of actions are not instantaneous in general \cite{delay1,delay2,delay3,delay4}. It may take a certain time delay. This leads to delayed and trend-based stochastic differential equations of mean-field type. 
\begin{equation} \label{delay:problem}
\begin{array}{ll} 
R(x,u) = -q_T var(x_T)  + \frac{(\bar{x}_T+\bar{\eta} \bar{y}_T)^{\rho}}{\rho}\\  %-r var(u_1) +
\int_0^T -q \ var(x) -r_1\ var(u_1) + \bar{r}_1  \frac{\bar{u}_1^{\rho}}{\rho}  dt,\\
\sup_{(u_1,u_2)} \mathbb{E}[ R(x,u) ],\\
\mbox{subject to }\\
dx= [ -(u_1-\bar{u}_1)+b_1(x-\bar{x})+b_2\epsilon (x-\bar{x})(u_2-\bar{u}_2)] dt \\
+[ -\bar{u}_1+\bar{b}_{11}\bar{x}+ \bar{b}_{12}\bar{y}+\bar{b}_{13}\bar{z} +\bar{b}_2\bar{u}_2 \bar{x} ] dt \\
+ \sigma(x-\bar{x})(u_2-\bar{u}_2)dB+ \bar{\sigma}\bar{x}\bar{u}_2 dB_o,\\
x(t')=x_0(t'), \ t'\in (-\delta, 0],\ s(0)=s_0,
\end{array}
\end{equation}
where $var(X)$ denotes the variance of the random variable $X,$ and  $\bar{X}(t)=\mathbb{E}[ X(t)| \ \mathcal{F}^{s,B_o}]$  is the conditional expectation with respect to the common noises $s,B_o.$ 

\begin{lemma}
The conditional expected trend $\bar{y}$ satisfies the following stochastic differential equation:
$$d\bar{y}=  [\bar{x}-\tau \bar{y}-e^{-\lambda \tau} \bar{z}]dt.$$ \eod
\end{lemma}

{\bf Proof: }

$$
dy= dt \int_{-\tau}^{0} e^{\lambda s}[\frac{d}{dt}x(t+s)] ds
%%= [e^{\lambda s} x(t+s)]-  \delta \int_{-\delta}^{0} e^{\lambda s}x(t+s)ds
= [x(t)-e^{-\lambda \tau} z(t) -\tau y(t)]dt.
$$
Taking the conditional expected values one obtains
$$d\bar{y}=  [\bar{x}-\tau \bar{y}-e^{-\lambda \tau} \bar{z}]dt.$$ This completes the proof.   \eod

\begin{prop} \label{semi:prop8}
The equilibrium strategies and the equilibrium payoff of the delayed MFTG (\ref{delay:problem}) are given by 
\begin{equation} \label{eq:delay12}
\begin{array}{ll} 
u^*_1= \sum_{s\in\mathcal{S}} \ind_{\{s(t)=s \}} [ \frac{\alpha}{r_1}(x-\bar{x})+(\bar{x}+\bar{\eta} \bar{y}) (\frac{\beta}{\bar{r}_1})^{\frac{1}{\rho-1}}], \\
u^*_2=\sum_{s\in\mathcal{S}} \ind_{\{s(t)=s \}} 
[ -\frac{b_2\epsilon}{\sigma^2}-\frac{(\bar{x}+\bar{\eta} \bar{y})(\beta\bar{b}_2+\beta_{B_o}\bar{\sigma})}{2(\rho-1)\beta\bar{\sigma}^2\bar{x}}], \\

\mathbb{E}[ R(x,u^*) ]=\mathbb{E}[ -\alpha(0,s_0) var(x)  + \beta(0,s_0) \frac{(\bar{x}_0+\bar{\eta} \bar{y}_0)^{\rho}}{\rho}],
 \\
\end{array}  
\end{equation}   
whenever the following system: 
\begin{equation} \label{eq:delay13}
\begin{array}{ll} 
d\alpha+
(2b_1\alpha+q+\sum_{s'}[\alpha(t,s')-\alpha(t,s)]\tilde{q}_{ss'}-\frac{\alpha^2}{r_1}-b_2^2\epsilon\alpha)dt=0,\\
\alpha(t,s)= q_T,\\
d\beta -(\frac{\rho}{\rho-1} \frac{(\beta\bar{b}_2+\beta_{B_o}\bar{\sigma})^2}{4\beta\bar{\sigma}^2} )dt\\
+\beta\rho(b_{11}+b_{13}e^{\lambda \tau})
+ \beta^{\frac{\rho}{\rho-1}}(\rho-1) (\bar{r}_1)^{-\frac{1}{\rho-1}} dt\\
+ \sum_{s'}[\beta(t,s')-\beta(t,s)]\tilde{q}_{ss'} dt\\
%-\frac{\rho}{\rho-1} \frac{(\beta\bar{b}_2+\beta_{B_o}\bar{\sigma})}{2\bar{\sigma}} dB_o
=0,\\
\beta(T,s)=1,\\
\bar{\eta}:=b_{13}e^{\lambda \tau}, 
\end{array} 
\end{equation}
has a unique solution $\alpha,\beta, \beta_{B_o}$ \eod
\end{prop}

 Note that the system in $\alpha$ has a positive solution if $q\geq 0, q_T\geq 0,\ r_1>0.$
With single regime $\mathcal{S}=\{s_0\}$ the $\beta$  equation yields

\begin{equation} \label{eq:delay14}
\begin{array}{ll} 
\dot{\beta}+ \beta[- \frac{\rho}{\rho-1} \frac{\bar{b}_2^2}{4\bar{\sigma}^2} 
+\rho(b_{11}+b_{13}e^{\lambda \tau})]
+ \beta^{\frac{\rho}{\rho-1}}(\rho-1) (\bar{r}_1)^{-\frac{1}{\rho-1}}=0,\\
\beta(T,s)=1,
\end{array} 
\end{equation}  
This is completely solvable with an  explicit solution given by
$$
\beta(t,s)=\mathbb{E}[   ( (1-\frac{c}{\omega}) e^{\frac{-\omega}{\rho-1}(T-t)}  +  \frac{c}{\omega}           )^{1-\rho} | \ s(t)=s]
$$
where  $$c:=(1-\rho)\bar{r}_1^{\frac{1}{1-\rho}},\ 
\omega:= \frac{\rho}{1-\rho} \frac{\bar{b}_2^2}{4\bar{\sigma}^2} 
+\rho(b_{11}+b_{13}e^{\lambda \tau})>0.$$
{\color{black}

{\bf Proof: \ }

Let $ f(x)= -\alpha\ var(x)  + \beta \frac{(\bar{x}+\bar{\eta} \bar{y})^{\rho}}{\rho},$ be a guess functional. 

\begin{equation} \label{eq:delay}
\begin{array}{ll} 
\mathbb{E}[ R-f(0)]= 
 \mathbb{E}-(q_T-\alpha_T)var(x_T)+(1-\beta_T)\frac{(\bar{x}+\bar{\eta} \bar{y})^{\rho}}{\rho}\\
+ \mathbb{E}\int_0^T -\left\{ d\alpha+(2b_1\alpha+q)dt   \right.\\ \left.
+\sum_{s'}[\alpha(t,s')-\alpha(t,s)]\tilde{q}_{ss'} dt+(-\frac{\alpha^2}{r_1}-b_2^2\epsilon\alpha) dt \right\}(x-\bar{x})^2 \\
-r_1[ u_1-\bar{u}_1 -\frac{\alpha}{r}(x-\bar{x}) ]^2 dt\\
 -\alpha\sigma^2[u_2-\bar{u}_2+\frac{b_2\epsilon}{\sigma^2}]^2 (x-\bar{x})^2 dt\\
+ \{   d\beta- \frac{\rho}{\rho-1} \frac{(\beta\bar{b}_2+\beta_{B_o}\bar{\sigma})^2}{4\beta\bar{\sigma}^2}    dt
-\frac{\rho}{\rho-1} \frac{(\beta\bar{b}_2+\beta_{B_o}\bar{\sigma})}{2\bar{\sigma}} dB_o\\
+(\beta\rho(b_{11}+b_{13}e^{\lambda \tau})
+ \beta^{\frac{\rho}{\rho-1}}(\rho-1) (\bar{r}_1)^{-\frac{1}{\rho-1}} )dt\\
+ \sum_{s'}[\beta(t,s')-\beta(t,s)]\tilde{q}_{ss'} dt  \} \frac{(\bar{x}+\bar{\eta} \bar{y})^{\rho}}{\rho} \\
+(\rho-1)\beta\bar{\sigma}^2\bar{x}^2 (\bar{x}+\bar{\eta} \bar{y})^{\rho-2} [\bar{u}_2+ \frac{(\bar{x}+\bar{\eta} \bar{y})(\beta\bar{b}_2+\beta_{B_o}\bar{\sigma})}{2(\rho-1)\beta\bar{\sigma}^2\bar{x}}  ]^2 dt\\
+[-(\bar{x}+\bar{\eta} \bar{y})^{\rho-1}\beta \bar{u}_1+\frac{\bar{r}_1 \bar{u}^{\rho}_1}{\rho}  \\
-\beta^{\frac{\rho}{\rho-1}} (\bar{x}+\bar{\eta} \bar{y})^{\rho}(1-\frac{1}{\rho}) (\bar{r}_1)^{-\frac{1}{\rho-1}}] dt,
\end{array}
\end{equation}
with the following careful matching $\eta=\bar{b}_{13} e^{\lambda \tau},\  \bar{b}_{12} = \bar{b}_{13}e^{\lambda \tau}(b_{11}+\lambda+\bar{b}_{13}e^{\lambda \tau}).$
The joint optimization over $(u_1,u_2)$ together with the mean-field terms $(\bar{u}_1,\bar{u}_2)$ gives the announced result provided that $\rho<1.$

\subsection{Mean-Field of MFTG}
This subsection we examine a class of mean-field of mean-field-type games. 

In view of the delayed mean-field-type game (\ref{delay:problem}), we have modified  $\bar{r}_1$  to be $\bar{r}_1(m)$ where $m$ is the conditional total consumption of the large population. Then, $m$ is obtained as
$$
m=\int \bar{u}_1 (t,\bar{x}, \bar{y}) \mu^{m}(t,d\bar{x}, d\bar{y}),$$
where $\mu^m(t,d\bar{x}, d\bar{y})$ is the conditional distribution of all players' states and  trends in the large population under $m,$
which  reduces to the fixed-point problem
$m= (\bar{x}(m)+\bar{\eta} \bar{y}(m)) (\frac{\beta(m)}{\bar{r}_1(m)})^{\frac{1}{\rho-1}}.$

Now consider the  following modified Cournot-Ross game with $I$ producers and a large population of potential consumers. The mean-field-type version of the game under common noise is analyzed in  \cite{hongkong}.
\begin{equation} \label{mfdelay:problem}
\begin{array}{ll} 
R_i(x,u) = \bar{R}-\frac{(x_T-\bar{x}_T)^{2k}}{2k}  + \frac{(\bar{x}_T)^{2k}}{2k}\\ +
\int_0^T - \bar{r}_i  \frac{(u_i-\bar{u}_i)^{2k}}{2k}+   {\color{blue}  \bar{x}^{2k-1} \bar{u}_i-\bar{r}_i  \frac{\bar{u}_i^{2k}}{2k}}  dt,\\
\sup_{u_i} \mathbb{E}[ R_i(x,u) ],\\
\mbox{subject to }\\
dx= [ D(m)-S] dt + \sigma(x-\bar{x})dB,\\
x(0)=x_0, \ s(0)=s_0,
\end{array}
\end{equation}
Let  $D(m)$ be the demand  generated by a large population of consumers. Given a demand $D(m),$ each macro-player $i$ has a certain utility of mean-field type. The
payoff function in the Cournot game {\color{blue} $\bar{x} \bar{u}_i+\bar{r}_i  \frac{\bar{u}_i^{2}}{2}$ } is modified to 
be {\color{blue} $  \bar{x}^{2k-1} \bar{u}_i-\bar{r}_i  \frac{\bar{u}_i^{2k}}{2k}$ } and some extra mean-field dependent terms. 

 By means of a direct method one can fully characterize the  mean-field equilibrium of (\ref{mfdelay:problem}). It is given by following set of equations:
\begin{equation} \label{mfdelay:solution}
\begin{array}{ll} 
\mathbb{E}[R_i(x,u^*)] =\mathbb{E}[ \bar{R}-\alpha_i(0,s_0)\frac{(x_0-\bar{x}_0)^{2k}}{2k}  - \bar{\alpha}_i(0,s_0)\frac{(\bar{x}_T)^{2k}}{2k}],\\
u^*_i =          (\frac{\alpha_i}{r_i})^{\frac{1}{2k-1}} (x-\bar{x} )        
  +    (\frac{\bar{\alpha}_i+1}{\bar{r}_i})^{\frac{1}{2k-1}} \bar{x}, \\
D(m)= [ \sum_{i}(\frac{\bar{\alpha}_i+1}{\bar{r}_i})^{\frac{1}{2k-1}}] \bar{x}.
\end{array}
\end{equation}
where $\alpha_i, \bar{\alpha}_i$ solve  a system of ordinary differential equations.

%\begin{remark}
%Note that in the mean-field game literature it is often considered an infinite number producers and  the supply process is scaled (averaged) $\frac{1}{n}\sum_{i}u_i. $
%While it seems reasonable to consider relatively large number of consumers for certain goods, 
%the assumption of infinite number of true decision-makers who are suppliers is not observed in practice. Another difficulty with the scaling model is that in the infinite regime the average supply becomes 
%deterministic which does not seem to capture stochastic demand observed in the markets. 
%
%Only few suppliers do have an effect on the price. In this example,  we have considered
% finite number of suppliers and arbitrary large number of consumers. We do not scale the supply.  The supply quantity $S=\sum_{i}u_i$ has to be matched with the total demand of the consumers $D(m).$ 
%\end{remark}

\section{MFTG beyond Brownian motions and Poisson }  \label{sec:mf:gv}
In this section class of mean-field-type games with a state dependent Gauss-Volterra noise is formulated and solved with a polynomial and mean-field dependent
 payoff for an arbitrary number of players and a finite time horizon. The control strategies are  linear state and mean-field feedbacks. 
A mean-field-type Nash equilibrium is verified for the game and the  optimal strategies are obtained using a direct method that does not require solving nonlinear partial integro-differential equations or forward-backward stochastic differential equations.  The example below is inspired from  \cite{refpp11,refpp12,refpp13,refpp14,refpp15,refpp16}. We add mean-field terms to these previous works. This will allow  us to solve variance or higher moment reduction problems.
 
 \subsection{Noncooperative MFTG under Gauss-Volterra processes}
This section examines a class of noncooperative mean-field-type games with non-quadratic cost and state driven by Gauss-Volterra processes.
\begin{equation} \label{gv:problem}
\begin{array}{ll} 
L_i(x,u) = q_{iT}\frac{(x_T-\bar{x}_T)^{2k_i}}{2k_i}  + \bar{q}_{iT} \frac{\bar{x}_T^{2\bar{k}_i}}{2\bar{k}_i}\\ +
\int_0^T q_{i}\frac{(x-\bar{x})^{2k_i}}{2k_i} +{r}_i  \frac{(u_i-\bar{u}_i)^{2{k}_i}}{2{k}_i}
+     \bar{q}_{i}\frac{\bar{x}^{2\bar{k}_i}}{2\bar{k}_i} + \bar{r}_i  \frac{\bar{u}_i^{2\bar{k}_i}}{2\bar{k}_i}  dt,\\
\inf_{u_i} \mathbb{E}[ L_i(x,u) ],\\
\mbox{subject to }\\
dx= [ b_1 (x-\bar{x})+\sum_{j}b_{2j} (u_j-\bar{u}_j)  +\bar{b}_1 \bar{x}+\sum_{j}\bar{b}_{2j} \bar{u}_j] dt\\
+(x-\bar{x}) [ \sigma dB+\int_{\Theta} \mu d\tilde{N}+\sigma_{gv} {\color{blue}dB_{gv}}],\\
\mathbb{P}(s(t+\epsilon)=s' | s,u)=\int_t^{t+\epsilon} \tilde{q}_{ss'} dt' + o(\epsilon),\ s'\neq s\\
x(0)=x_0, \ s(0)=s_0,
\end{array}
\end{equation}
where $k_i\geq 1 ,\bar{k}_i\geq 1$ are natural numbers, the coefficients are time and switching dependent, 

\begin{remark}The cost functional is clearly non-quadratic for $ k_i> 1 $ or $\bar{k}_i>1.$  
For $k_i=1$ the equilibrium of the variance reduction game (\ref{gv:problem}) under Gauss-Volterra processes is obtained.
\end{remark}

\begin{prop} \label{semi:prop9} Assume $q_i>0, \bar{q}_i>0, r_i, \bar{r}_i>\delta$  and $\int_{\theta}[ (1+\mu)^{2k_i}-1 -2k_i\mu   ] \nu(d\theta) <+\infty.$
The mean-field Nash equilibrium of the mean-field type game (\ref{gv:problem}) under Gauss-Volterra process is 
given by 
\begin{equation}
\begin{array}{ll}
u^*_i= \sum_{s\in\mathcal{S}} \ind_{\{s(t)=s \}} 
[\Bigg( -\frac{b_{2i}\alpha_i}{r_i}\Bigg)^{\frac{1}{2k_i-1}} (x-\bar{x})+\Bigg( -\frac{\bar{b}_{2i}\bar{\alpha}_i}{\bar{r}_i}\Bigg)^{\frac{1}{2\bar{k}_i-1}} \bar{x}],\ \\
\mathbb{E}L_i=\mathbb{E}[ \alpha_i(0,s_0) \frac{(x_0-\bar{x}_0)^{2k_i}}{2k_i}  + \bar{\alpha}_i (0,s_0)\frac{\bar{x}_0^{2\bar{k}_i}}{2\bar{k}_i}],
\end{array}
\end{equation}
whenever the following system of ordinary differential equations admit a positive solution which does not blowup within 
$[0,T].$

\begin{equation} \label{fieq8}
\begin{array}{ll}
\dot{\alpha}_i  + q_{i}+ 2k_i\alpha_ib_1 + \alpha_i k_i(2k_i-1)(\sigma^2+{\color{blue} \sigma^2_{cogv}} )\\
+ {\alpha}_i\int_{\theta}[ (1+\mu)^{2k_i}-1 -2k_i\mu   ] \nu(d\theta)\\
+(\sum_{s'} [{\alpha}_i(t,s')- \alpha_i(t,s)]\tilde{q}_{ss'})  \\ 
-(2k_i-1)  r_i ( -\frac{b_{2i}\alpha_i}{r_i})^{\frac{2k_i}{2k_i-1}}\\
+2k_i\alpha_i [\sum_{j\neq i} b_{2j} ( -\frac{b_{2j}\alpha_j}{r_j})^{\frac{1}{2k_j-1}}=0,\\
{\alpha}_i(T,s)=q_i(T,s),\\
 \dot{\bar{\alpha}}_i + \bar{q}_{i} + 2\bar{k}_i\bar{\alpha}_i \bar{b}_1  
+(\sum_{s'} [\bar{\alpha}_i(t,s')- \bar{\alpha}_i(t,s)]\tilde{q}_{ss'})\\ 
- (2\bar{k}_i-1) \bar{r}_i ( -\frac{\bar{b}_{2i}\bar{\alpha}_i}{\bar{r}_i})^{\frac{2\bar{k}_i}{2\bar{k}_i-1}} \\
+2\bar{k}_i\bar{\alpha}_i[\sum_{j\neq i}  \bar{b}_{2j} ( -\frac{\bar{b}_{2j}\bar{\alpha}_j}{\bar{r}_j})^{\frac{1}{2\bar{k}_j-1}} ] 
=0,\\
\bar{\alpha}_i(T,s)=\bar{q}_i(T,s),\\

%%\sigma_{cogv}^2=\frac{d}{dt}\bigg[\int_0^t   \Bigg\{ K(t'_+,t')\sigma_{gv}(t')+\int_{t'}^{t}\sigma_{gv}(t'')K(t',t'')dt''\Bigg\}^2  dt'\bigg],

\end{array}
\end{equation}
\eod
\end{prop}

{\bf Proof: }

Consider the guess functional $f_i=\alpha_i \frac{(x-\bar{x})^{2k_i}}{2k_i}  + \bar{\alpha}_i \frac{\bar{x}^{2\bar{k}_i}}{2\bar{k}_i}$
\begin{equation} \label{gv:expand}
\begin{array}{ll} 
\mathbb{E}[L_i- f_i(0)]= \mathbb{E}( {q}_{iT} -\alpha_i(T,s(T))) \frac{(x_T-\bar{x}_T)^{2k_i}}{2k_i} \\ 
+( \bar{q}_{iT} - \bar{\alpha}_i(T,s(T))) \frac{\bar{x}_T^{2\bar{k}_i}}{2\bar{k}_i}\\
+
\mathbb{E}\int_0^T q_{i}\frac{(x-\bar{x})^{2k_i}}{2k_i} 
\dot{\alpha}_i \frac{(x-\bar{x})^{2k_i}}{2k_i} 
 + \alpha_ib_1 (x-\bar{x})^{2k_i}\\
+{r}_i  \frac{(u_i-\bar{u}_i)^{2\bar{k}_i}}{2\bar{k}_i}+[\sum_{j} \alpha_i (x-\bar{x})^{2k_i-1}b_{2j} (u_j-\bar{u}_j)] \\ 
+\alpha_i k_i(2k_i-1)(\sigma^2{\color{blue}+\sigma_{cogv}^2})  \frac{(x-\bar{x})^{2k_i}}{2k_i}\\
+ {\alpha}_i\int_{\theta}[ (1+\mu)^{2k_i}-1 -\mu 2k_i  ] \nu(d\theta) \frac{(x-\bar{x})^{2k_i}}{2k_i}\\
+(\sum_{s'} [{\alpha}_i(t,s')- \alpha_i(t,s)]\tilde{q}_{ss'}) \frac{(x-\bar{x})^{2k_i}}{2k_i} \\
%%\bar{\alpha}_i(T)) \frac{\bar{x}^{2\bar{k}}}{2\bar{k}}
+     \bar{q}_{i}\frac{\bar{x}^{2\bar{k}_i}}{2\bar{k}_i}  +
+\dot{\bar{\alpha}}_i \frac{\bar{x}^{2\bar{k}_i}}{2\bar{k}_i}
+\bar{\alpha}_i \bar{x}^{2\bar{k}_i} \bar{b}_1\\
+ \bar{r}_i  \frac{\bar{u}_i^{2\bar{k}_i}}{2\bar{k}_i}+\bar{\alpha}_i \bar{x}^{2\bar{k}_i-1} [\sum_{j}\bar{b}_{2j}\bar{u}_j]\\
+(\sum_{s'} [\bar{\alpha}_i(t,s')- \bar{\alpha}_i(t,s)]\tilde{q}_{ss'})  \frac{\bar{x}^{2\bar{k}_i}}{2\bar{k}_i}
\end{array}
\end{equation}

We complete the following term: 
\begin{equation} 
\begin{array}{ll} 
{r}_i  \frac{(u_i-\bar{u}_i)^{2\bar{k}}}{2{k}_i}+  \alpha_i (x-\bar{x})^{2k_i-1}b_{2i} (u_i-\bar{u}_i)\\
={r}_i \frac{(u_i-\bar{u}_i)^{2\bar{k}}}{2{k}_i}+   \alpha_i (x-\bar{x})^{2k_i-1}b_{2i} (u_i-\bar{u}_i)\\
= {r}_i \frac{(u_i-\bar{u}_i)^{2\bar{k}_i}}{2{k}_i}+   \alpha_i (x-\bar{x})^{2k_i-1}b_{2i} (u_i-\bar{u}_i)\\
+(2k_i-1) r_i ( -\frac{b_{2i}\alpha_i}{r_i})^{\frac{2k_i}{2k_i-1}} \frac{(x-\bar{x})^{2k_i}}{2k_i}\\
- (2k_i-1) r_i ( -\frac{b_{2i}\alpha_i}{r_i})^{\frac{2k_i}{2k_i-1}} \frac{(x-\bar{x})^{2k_i}}{2k_i}\\
\end{array}
\end{equation} 
A similar completion is done for the terms  in $\bar{u}_i.$
Thus,
\begin{equation} \label{gv:expandtyp}
\begin{array}{ll} 
\mathbb{E}[L_i- f_i(0)]\\ = \mathbb{E}( {q}_{iT} -\alpha_i(T)) \frac{(x_T-\bar{x}_T)^{2k_i}}{2k_i}  
+( \bar{q}_{iT} - \bar{\alpha}_i(T)) \frac{\bar{x}_T^{2\bar{k}_i}}{2\bar{k}_i}\\

+\mathbb{E}\int_0^T
\left\{ [\dot{\alpha}_i  + q_{i}+ 2k_i\alpha_ib_1 + \alpha_i k_i(2k_i-1)(\sigma^2+\sigma^2_{cogv} ) \right.\\
+ {\alpha}_i\int_{\theta}[ (1+\mu)^{2k_i}-1 -\mu 2k_i  ] \nu(d\theta)\\
+\sum_{s'} [{\alpha}_i(t,s')- \alpha_i(t,s)]\tilde{q}_{ss'}  \\ 
-(2k_i-1)  r_i ( -\frac{b_{2i}\alpha_i}{r_i})^{\frac{2k_i}{2k_i-1}}\\ +
2k_i\alpha_i \sum_{j\neq i} b_{2j} ( -\frac{b_{2j}\alpha_j}{r_j})^{\frac{1}{2k_j-1}} \left. \right\} \frac{(x-\bar{x})^{2k_i}}{2k_i} 
\\
+  \left\{   \dot{\bar{\alpha}}_i + \bar{q}_{i} + 2\bar{k}_i\bar{\alpha}_i \bar{b}_1  
+(\sum_{s'} [\bar{\alpha}_i(t,s')- \bar{\alpha}_i(t,s)]\tilde{q}_{ss'})\right. \\
- (2\bar{k}_i-1) \bar{r}_i ( -\frac{\bar{b}_{2i}\bar{\alpha}_i}{\bar{r}_i})^{\frac{2\bar{k}_i}{2\bar{k}_i-1}} \\
+2\bar{k}_i\bar{\alpha}_i[\sum_{j\neq i}  \bar{b}_{2j} ( -\frac{\bar{b}_{2j}\bar{\alpha}_j}{\bar{r}_j})^{\frac{1}{2\bar{k}_j-1}} ]  \left. \right\} \frac{\bar{x}^{2\bar{k}_i}}{2\bar{k}_i}\\ 
+{r}_i [\frac{(u_i-\bar{u}_i)^{2{k}_i}}{2{k}_i}+   \alpha_i (x-\bar{x})^{2k_i-1}b_{2i} (u_i-\bar{u}_i)]\\
+ (2k_i-1) r_i ( -\frac{b_{2i}\alpha_i}{r_i})^{\frac{2k_i}{2k_i-1}} \frac{(x-\bar{x})^{2k_i}}{2k_i}\\
+\bar{r}_i [\frac{\bar{u}_i^{2\bar{k}_i}}{2\bar{k}_i}+   \bar{\alpha}_i \bar{x}^{2\bar{k}_i-1}b_{2i} \bar{u}_i]\\
+(2\bar{k}_i-1) \bar{r}_i ( -\frac{\bar{b}_{2i}\bar{\alpha}_i}{\bar{r}_i})^{\frac{2\bar{k}}{2\bar{k}_i-1}} \frac{\bar{x}^{2\bar{k}_i}}{2\bar{k}_i}\\
\end{array}
\end{equation}   

Noting  for $k_i>\frac{1}{2}, r_i>0$ one has:
\begin{equation} \label{equiif}
\begin{array}{ll}
{r}_i \frac{(u_i-\bar{u}_i)^{2\bar{k}_i}}{2\bar{k}_i}+   \alpha_i (x-\bar{x})^{2k_i-1}b_{2i} (u_i-\bar{u}_i)\\
+ (2k_i-1) r_i ( -\frac{b_{2i}\alpha_i}{r_i})^{\frac{2k_i}{2k_i-1}} \frac{(x-\bar{x})^{2k_i}}{2k_i}\geq 0,\\
\bar{r}_i \frac{\bar{u}_i^{2\bar{k}}}{2\bar{k}_i}+   \bar{\alpha}_i \bar{x}^{2\bar{k}_i-1}b_{2i} \bar{u}_i\\
+ (2\bar{k}_i-1) \bar{r}_i ( -\frac{\bar{b}_{2i}\bar{\alpha}_i}{\bar{r}_i})^{\frac{2\bar{k}_i}{2\bar{k}_i-1}} \frac{\bar{x}^{2\bar{k}_i}}{2\bar{k}_i}
\geq 0,
\end{array}
\end{equation} 
 with equalities in (\ref{equiif}) iff 
 $$
 u_i-\bar{u}_i= ( -\frac{b_{2i}\alpha_i}{r_i})^{\frac{1}{2k_i-1}} (x-\bar{x}),\ \  \bar{u}_i=
 ( -\frac{\bar{b}_{2i}\bar{\alpha}_i}{\bar{r}_i})^{\frac{1}{2\bar{k}_i-1}} \bar{x}.
 $$
 By identification, one obtains the announced result.  \eod

 The derived system of equations (\ref{fieq8}) are inhomogeneous differential system where it is known that existence and uniqueness, nonexistence or nonuniqueness may occur.

\subsubsection*{Existence}
Some results on sufficient conditions for
the existence of trajectories satisfying the associated set of non-linear differential equations  (\ref{fieq8}) are
outlined. Below we present Carath\'eodory  conditions for existence of a solution.

Here, the non-linear differential system  (\ref{fieq8}) can be written as 
\begin{equation} \label{caratheodory}
\begin{array}{ll}
\dot{\alpha}= h(t,\alpha),\ \ \\
\alpha_i(T,s)=q_i(T,s),\ s\in\mathcal{S}, i\in \mathcal{N}.
\end{array}
\end{equation} 
Assume that 
\begin{itemize}
\item For each fixed time $t\in [0,T],$   $h$ is continuous in $\alpha.$
\item For each fixed $\alpha,$   $h$ is measurable in $t.$
\item Given a nonempty compact $C\subset  \mathbb{R}^{n |\mathcal{S}|}$ and interval $[T-\epsilon,T],$ there is an integrable positive function $\hat{h}$ on the time interval $[T-\epsilon,T]$ such that
$|h(t,\alpha)|\leq \hat{h}(t),\  $ for all $(t,\alpha)\in [0,T]\times C.$
%$$|h(t,\alpha)\leq \hat{h}(t),\  |h_{\alpha}(t,\alpha)\leq \hat{h}(t), $$ 
\end{itemize}

The interval of definition of the solution depends on the terminal value $q(T,s).$ For each terminal condition  $q(T,.)\in \mathbb{R}^{n |\mathcal{S}|},$  there is an interval $(T-\epsilon, T), \ \epsilon>0$ where this  non-linear differential system (\ref{caratheodory}) has at least one  solution. We refer to    \cite[Theorem 1.1, Chapter 2]{cod} for a detailed proof.

\subsubsection*{Uniqueness}
It is a well-known
result that not every non-linear differential system has a unique solution. 
Therefore, the uniqueness issue is dealt with in a separate result. 
We provide two sufficient conditions for having at most one solution:
\begin{itemize}
\item  If   $h$ is continuously differentiable in $\alpha$ and $|h_{\alpha}(t,\alpha)|\leq \hat{h}(t), $  on $(t,\alpha)\in [0,T]\times C,$ then there is {\it at most one}  solution on $[T-\epsilon, T].$ 
\item  If $|h(t,\alpha_1)-h(t,\alpha_2)|\leq \hat{h}(t) |\alpha_1-\alpha_2|,$  on $(t,\alpha)\in [0,T]\times C,$ then there is {\it at most one}  solution on $[T-\epsilon, T].$ 
%\item
 % If $h$ is negatively monotone: $ \langle \alpha_1-\alpha_2, h(t,\alpha_1)-h(t,\alpha_2) \rangle
\end{itemize}

It is important the notice that the function $h$ is not necessarily globally Lipschitz in $\alpha.$ For example,  for $k\geq 1,$ $\alpha_i^{\frac{2k}{2k-1}}$ is not necessarily globally Lipschitz. Therefore we need estimates of $\alpha.$ We rely on the original dynamic optimization problem to derive  lower and upper  bounds on $\alpha_i.$  Since  $q(T,.), q(t,.)\succ 0, r_i>\delta>0$ by assumption, lower bound for $\alpha_i$ is zero. This can also be obtained directly from the problem formulation as the cost is positive.
By summing up (\ref{fieq8}) over $i\in\mathcal{I},$ an upper bound is obtained as  $\sum_{i}\alpha_i$ is bounded subject to integrability condition of the coefficients. 
 %A upper bound   

Note, however, that the stationary system may have multiple solutions, depending on the parameters. 
%%\subsubsection*{Existence and uniqueness}

\subsubsection*{Admissibility of the coefficient solution}
%Also, the existence result  does not say anything about the positiveness of the solution $\alpha.$  If $q(T,.)\succ 0,$ a lower bound is zero.

As $\epsilon$ in the Carath\'eodory existence result  depends on $q(T,.),$ the maximal interval in which the solution is defined may depend on $q(T,.).$
Thus, we need to examine the singularity of $\alpha$ in (\ref{fieq8}). In order for the control strategies to be admissible, we seek  sufficient conditions for non-blow-up (no escape) within $[0,T$
If  $q_i>0, r_i>\delta>0$  and all coefficients continuous, then there is no escape within  $[T-\epsilon,T].$  If in addition the coefficient functions $ b_1, \sigma^2, \sigma^2_{cogv},
\int_{\theta}[ (1+\mu)^{2k_i}-1 -2k_i\mu   ] \nu(d\theta),$   $r_i ( \frac{b_{2i}}{r_i})^{\frac{2k_i}{2k_i-1}} $ and
$(\frac{{b}_{2i}}{{r}_i})^{\frac{1}{2{k}_j-1}}$ are all integrable within $\mathcal{T}$, then there is no escape of $\alpha$ within the entire $[0,T]$ as the estimates of $\sum_{i\in\mathcal{I}}\alpha_i$ is finite in $\mathcal{T}.$

Similar reasoning works for $\bar{\alpha}$ when  $\bar{q}_i>0, \bar{r}_i>\delta>0$ and the coefficient functions
$ \bar{b}_1, $    $\bar{r}_i ( \frac{\bar{b}_{2i}}{\bar{r}_i})^{\frac{2\bar{k}_i}{2\bar{k}_i-1}},$ 
$(\frac{\bar{b}_{2i}}{r_i})^{\frac{1}{2\bar{k}_i-1}}$ are integrable within $\mathcal{T}.$

{ \color{blue}

At equilibrium, the mean-field term $\bar{x}$ in Proposition \ref{semi:prop9} solves 
\begin{equation}
\begin{array}{ll}
d\bar{x}= [ \bar{b}_1 +\sum_{j}\bar{b}_{2j} 
\Bigg( -\frac{\bar{b}_{2j}\bar{\alpha}_j}{\bar{r}_j}\Bigg)^{\frac{1}{2\bar{k}_j-1}} ]  \bar{x}dt,\\
\bar{x}(0)=\bar{x}_0,
\end{array}
\end{equation}
which admits a unique solution within $[0,T]$ subject to the integrability of the regime switching 
dependent coefficient
$[ \bar{b}_1 +\sum_{j}\bar{b}_{2j} 
\Bigg( -\frac{\bar{b}_{2j}\bar{\alpha}_j}{\bar{r}_j}\Bigg)^{\frac{1}{2\bar{k}_j-1}} ]$ over $[0,T].$
}
\subsection{Fully Cooperative MFTG under Gauss-Volterra Noise}
In this subsection we choose $k_i=k,\ \bar{k}_i=\bar{k}$  and assume that the $I$ decision-makers are fully cooperative. They jointly decide and solve the following problem:

\begin{equation} \label{gv:problemcoop}
\begin{array}{ll} 
\inf_{(u_1,\ldots, u_n)} \mathbb{E}[ \sum_{i\in \mathcal{I}} L_i(x,u) ],\\
\mbox{subject to }\\
dx= [ b_1 (x-\bar{x})+\sum_{j}b_{2j} (u_j-\bar{u}_j)  +\bar{b}_1 \bar{x}+\sum_{j}\bar{b}_{2j} \bar{u}_j] dt\\
+(x-\bar{x}) [ \sigma dB+\int_{\Theta} \mu d\tilde{N}+\sigma_{gv} {\color{blue}dB_{gv}}],\\
\mathbb{P}(s(t+\epsilon)=s' | s,u)=\int_t^{t+\epsilon} \tilde{q}_{ss'} dt' + o(\epsilon),\ s'\neq s\\
x(0)=x_0, \ s(0)=s_0,
\end{array}
\end{equation}

\begin{prop}  \label{globalpropyy}
The global optimum of the fully cooperative mean-field type game (\ref{gv:problemcoop}) under Gauss-Volterra process is 
given by 
\begin{equation}
\begin{array}{ll}
u^*_i= \sum_{s\in\mathcal{S}} \ind_{\{s(t)=s \}} 
[\Bigg( -\frac{b_{2i}\alpha_0}{r_i}\Bigg)^{\frac{1}{2k-1}} (x-\bar{x})+\Bigg( -\frac{\bar{b}_{2i}\bar{\alpha}_0}{\bar{r}_i}\Bigg)^{\frac{1}{2\bar{k}-1}} \bar{x}],\ \\
\mathbb{E}L_i=\mathbb{E}[ \alpha_0(0,s_0)\frac{(x_0-\bar{x}_0)^{2k}}{2k}  + \bar{\alpha}_0 (0,s_0) \frac{\bar{x}_0^{2\bar{k}}}{2\bar{k}}],
\end{array}
\end{equation}
whenever the following system of ordinary differential equations admit a positive solution which does not blowup within 
$[0,T].$

\begin{equation}
\begin{array}{ll}
\dot{\alpha}_0  + (\sum_{i}q_{i})+ 2k\alpha_0b_1 + \alpha_0 k(2k-1)(\sigma^2{\color{blue}+\sigma^2_{cogv} })\\
+ {\alpha}_0\int_{\theta}[ (1+\mu)^{2k}-1 -2k\mu   ] \nu(d\theta)\\
+\sum_{s'} [{\alpha}_0(t,s')- \alpha_0(t,s)]\tilde{q}_{ss'}  \\
-(2k-1)  \sum_{i}r_i ( -\frac{b_{2i}\alpha_0}{r_i})^{\frac{2k}{2k-1}}
=0,\\
{\alpha}_0(T,s)=\sum_i q_i(T,s),\\
 \dot{\bar{\alpha}}_0+ (\sum_{i}\bar{q}_{i}) + 2\bar{k}\bar{\alpha}_0 \bar{b}_1  
+(\sum_{s'} [\bar{\alpha}_0(t,s')- \bar{\alpha}_0(t,s)]\tilde{q}_{ss'}) \\
- (2\bar{k}-1) \sum_{i}\bar{r}_i ( -\frac{\bar{b}_{2i}\bar{\alpha}_0}{\bar{r}_i})^{\frac{2\bar{k}}{2\bar{k}-1}} 
=0,\\
\bar{\alpha}_0(T,s)=\sum_{i}\bar{q}_i(T,s),\\
\end{array}
\end{equation} 
\eod
\end{prop}  

{\bf Proof: } The proof follows similar steps as in Proposition \ref{semi:prop9}.

\begin{remark} A sufficient condition for existence and uniqueness of the global optimum of mean-field type  is obtained for
 $r_i>\delta>0,\ \sum_{i}q_i>0, \int_{\theta}[ (1+\mu)^{2k}-1 -2k\mu   ] \nu(d\theta) < +\infty$  and all coefficients continuous. Then there is no escape within  $[T-\epsilon,T].$  If in addition the coefficient functions $ b_1, \sigma^2, \sigma^2_{cogv},
\int_{\theta}[ (1+\mu)^{2k}-1 -2k \mu   ] \nu(d\theta),$   $r_i ( \frac{b_{2i}}{r_i})^{\frac{2k}{2k-1}} $ and
$(\frac{{b}_{2i}}{{r}_i})^{\frac{1}{2{k}-1}}$ are all integrable within $[0,T]$, then is no escape of $\alpha$ within the entire interval $[0,T].$ 
Similar reasoning works for $\bar{\alpha}$ when  $ \bar{r}_i>\delta>0,  i \in \mathcal{I}, \ \sum_{j}\bar{q}_j>0,$ and the coefficient functions
$ \bar{b}_1, $    $\bar{r}_i ( \frac{\bar{b}_{2i}}{\bar{r}_i})^{\frac{2\bar{k}_i}{2\bar{k}-1}},$ 
$(\frac{\bar{b}_{2i}}{r_i})^{\frac{1}{2\bar{k}-1}}$ are integrable within $[0,T].$

{ \color{blue} 
At the global optimum, the mean-field term $\bar{x}$  in Proposition \ref{globalpropyy} solves 
\begin{equation}
\begin{array}{ll}
d\bar{x}= \left[ \bar{b}_1 +\sum_{j}\bar{b}_{2j} 
\Bigg( -\frac{\bar{b}_{2j}\bar{\alpha}_0}{\bar{r}_j}\Bigg)^{\frac{1}{2\bar{k}-1}} \right]  \bar{x}dt,\\
\bar{x}(0)=\bar{x}_0,
\end{array}
\end{equation}
which admits a unique solution within $[0,T]$ subject to the integrability of the regime switching 
dependent coefficient
$ \bar{b}_1 +\sum_{j}\bar{b}_{2j} 
\Bigg( -\frac{\bar{b}_{2j}\bar{\alpha}_0}{\bar{r}_j}\Bigg)^{\frac{1}{2\bar{k}-1}} $ over $[0,T].$
}
\end{remark}
{ \color{blue} 
\begin{remark}
Notice that the differential equation 
\begin{equation}
\begin{array}{ll}
\dot{\xi}+q+c_1\xi -c_2 \xi^{p}=0,\\
\xi(T)=q(T)\geq 0\\
 q(t)> 0, c_2(t)> 0\\
 p=\frac{2k}{2k-1}>1, 
\end{array}
\end{equation} has a unique solution within $[0,T].$ Moreover, the unique solution is positive.
\end{remark}
}
\subsection{Adversarial Mean-Field-Type Game under Gauss-Volterra Noise}
In this subsection we choose $k_i=k\geq 1,\ \bar{k}_i=\bar{k}\geq 1$  and assume that the $n$ decision-makers are 
divided into two teams $\mathcal{I}_{+}=\{ i\in \mathcal{I} \ | \   r_i>\delta>0, \bar{r}_i>\delta>0  \}, \mathcal{I}_{-}=\{ i\in \mathcal{I} | \  
 r_i<-\delta<0, \bar{r}_i<-\delta<0 \}$ and $\mathcal{I}=\mathcal{I}_{+}\cup \mathcal{I}_{-}.$
 The decision-makers in  team $\mathcal{I}_{+}$ minimize the functional
  $ \mathbb{E} L_{ad}(x,u)$ over $(u_i)_{i\in \mathcal{I}_{+}}.$
 The decision-makers in  team $\mathcal{I}_{-}$ maximize   $ \mathbb{E} L_{ad}(x,u)$ over $(u_j)_{j\in \mathcal{I}_{-}}$
This leads to a minmax game problem:

\begin{equation} \label{gv:problemad}
\begin{array}{ll} 
L_{ad}(x,u) = q_{T}\frac{(x_T-\bar{x}_T)^{2k}}{2k}  + \bar{q}_{T} \frac{\bar{x}_T^{2\bar{k}}}{2\bar{k}}\\ +
\int_0^T q\frac{(x-\bar{x})^{2k}}{2k} +\sum_{i}{r}_i  \frac{(u_i-\bar{u}_i)^{2\bar{k}}}{2\bar{k}}
+     \bar{q}\frac{\bar{x}^{2k}}{2k} + \sum_{i}\bar{r}_i  \frac{\bar{u}_i^{2\bar{k}}}{2\bar{k}}  dt,\\

\inf_{(u_i)_{i\in\mathcal{I}_+}} \sup_{(u_j)_{j\in\mathcal{I}_{-}}} \mathbb{E}[ L_{ad}(x,u) ],\\
\mbox{subject to }\\
dx= [ b_1 (x-\bar{x})+\sum_{j}b_{2j} (u_j-\bar{u}_j)  +\bar{b}_1 \bar{x}+\sum_{j}\bar{b}_{2j} \bar{u}_j] dt\\
+(x-\bar{x}) [ \sigma dB+\int_{\Theta} \mu d\tilde{N}+\sigma_{gv} {\color{blue}dB_{gv}}],\\
\mathbb{P}(s(t+\epsilon)=s' | s,u)=\int_t^{t+\epsilon} \tilde{q}_{ss'} dt' + o(\epsilon),\ s'\neq s\\
x(0)=x_0, \ s(0)=s_0,
\end{array}
\end{equation}

A mean-field-type risk-neutral saddle point is a strategy profile $(u^*_j, \ j\in I_{+}),$ of the team of defenders and  $(u^*_j, \ j\in I_{-})$ of the team of attackers such that 
\begin{align*}
\begin{array}{l}
\mathbb{E}L_{\mathrm{ad}}(x,s,(u^*_i)_{i\in I_{+} }, (u_j)_{j\in I_{-} })  \leq \mathbb{E} L_{\mathrm{ad}}(X,s,u^*) \\
\leq  \mathbb{E}L_{\mathrm{ad}}(X,s,(u_i)_{i\in I_{+} }, (u^*_j)_{j\in I_{-} }), \  \forall \ (u_i)_{i\in I_{+} },  (u_j)_{j\in I_{-} }
\end{array}
\end{align*}

\begin{prop} \label{addlpgv}   Assume $ \mathcal{I}_{+}\cup   \mathcal{I}_{-}= \mathcal{I}$ and $q>0 ,  \bar{q}> 0.$ Then,
the minmax solution of the adversarial mean-field type game (\ref{gv:problemad}) under Gauss-Volterra process is 
given by 
\begin{equation}
\begin{array}{ll}
i\in \mathcal{I}_{+},\\
u^*_i=\sum_{s\in\mathcal{S}} \ind_{\{s(t)=s \}}
[  \Bigg( -\frac{b_{2i}\alpha_{ad}}{r_i}\Bigg)^{\frac{1}{2k-1}} (x-\bar{x})+\Bigg( -\frac{\bar{b}_{2i}\bar{\alpha}_{ad}}{\bar{r}_i}\Bigg)^{\frac{1}{2\bar{k}-1}} \bar{x}],\ \\
j\in \mathcal{I}_{-},\\
u^*_j=\sum_{s\in\mathcal{S}} \ind_{\{s(t)=s \}} \Bigg( -\frac{b_{2j}\alpha_{ad}}{r_j}\Bigg)^{\frac{1}{2k-1}} (x-\bar{x})+\Bigg( -\frac{\bar{b}_{2j}\bar{\alpha}_{ad}}{\bar{r}_j}\Bigg)^{\frac{1}{2\bar{k}-1}} \bar{x},\ \\

\mathbb{E}L_{ad}=\mathbb{E}[ \alpha_{ad}(0,s_0)\frac{(x_0-\bar{x}_0)^{2k}}{2k}  + \bar{\alpha}_{ad} (0,s_0) \frac{\bar{x}_0^{2\bar{k}}}{2\bar{k}}],
\end{array}
\end{equation} 
whenever the following system of ordinary differential equations admit a positive solution which does not blowup within the horizon 
$[0,T].$

\begin{equation} 
\begin{array}{ll}
0=\dot{\alpha}_{ad}  + q+ 2k\alpha_{ad}b_1 + \alpha_{ad} k(2k-1)(\sigma^2 {\color{blue}+\sigma^2_{cogv}} )\\
+ {\alpha}_{ad}\int_{\theta}[ (1+\mu)^{2k}-1 -2k\mu   ] \nu(d\theta)\\
+\sum_{s'} [{\alpha}_{ad}(t,s')- \alpha_{ad}(t,s)]\tilde{q}_{ss'}  \\ 
{\color{blue} -(2k-1) \alpha_{ad}^{\frac{2k}{2k-1}} [  \sum_{i\in \mathcal{I}_{+}}r_i ( -\frac{b_{2i}}{r_i})^{\frac{2k}{2k-1}}
+ \sum_{j\in \mathcal{I}_{-}}r_j( -\frac{b_{2j}}{r_j})^{\frac{2k}{2k-1}} ]},\\
{\alpha}_{ad}(T,s)=q(T,s),\\
 0=\dot{\bar{\alpha}}_{ad}+ \bar{q} + 2\bar{k}\bar{\alpha}_{ad} \bar{b}_1  
+(\sum_{s'} [\bar{\alpha}_{ad}(t,s')- \bar{\alpha}_{ad}(t,s)]\tilde{q}_{ss'}) \\
{\color{blue}- (2\bar{k}-1)\bar{\alpha}_{ad}^{\frac{2k}{2k-1}} [
\sum_{i\in \mathcal{I}_{+}}\bar{r}_i ( -\frac{\bar{b}_{2i}}{\bar{r}_i})^{\frac{2\bar{k}}{2\bar{k}-1}} 
+\sum_{j\in \mathcal{I}_{-}}\bar{r}_j ( -\frac{\bar{b}_{2j}}{\bar{r}_j})^{\frac{2\bar{k}}{2\bar{k}-1}} 
]},\\
\bar{\alpha}_{ad}(T,s)=\bar{q}(T,s),\\
\end{array}
\end{equation}
In this case, the minmax solution is also a maxmin solution, hence $(u_i)_{i\in \mathcal{I}_{+}}, (u_j)_{j\in \mathcal{I}_{-}}$ is a saddle point. Thus, the adversarial mean-field-type game has a value $\mathbb{E}L_{ad}(x,(u_i)_{i\in \mathcal{I}_{+}}, (u_j)_{j\in \mathcal{I}_{-}}).$
\end{prop}
The proof follows similar steps as above by exploiting the strict convex-concave and coercivity properties of the cost functional.   

\begin{remark} A sufficient condition for existence and uniqueness of the minmax point of mean-field type  is obtained:
$$ 2k-1>0,\  q>0,
 \sum_{i\in \mathcal{I}_{+}}r_i ( \frac{b_{2i}}{r_i})^{\frac{2k}{2k-1}}
+ \sum_{j\in \mathcal{I}_{-}}r_j( \frac{b_{2j}}{r_j})^{\frac{2k}{2k-1}} >0,
$$
and  the coefficient functions $ b_1, \sigma^2, \sigma^2_{cogv},
\int_{\theta}[ (1+\mu)^{2k}-1 -2k \mu   ] \nu(d\theta),$   $r_i ( \frac{b_{2i}}{r_i})^{\frac{2k}{2k-1}} $ and
$(\frac{{b}_{2i}}{{r}_i})^{\frac{1}{2{k}-1}}$ are all integrable within $[0,T]$, then is no escape of $\alpha_{ad}$ within the entire interval $[0,T].$ 
 
Similar reasoning works for $\bar{\alpha}_{ad}$ when  $$2\bar{k}-1>0,\ \bar{q}>0, \sum_{i\in \mathcal{I}_{+}}\bar{r}_i ( \frac{\bar{b}_{2i}}{\bar{r}_i})^{\frac{2\bar{k}}{2\bar{k}-1}} 
+\sum_{j\in \mathcal{I}_{-}}\bar{r}_j ( \frac{\bar{b}_{2j}}{\bar{r}_j})^{\frac{2\bar{k}}{2\bar{k}-1}} >0,$$ 
and   the coefficient functions
$ \bar{b}_1, $    $\bar{r}_i ( \frac{\bar{b}_{2i}}{\bar{r}_i})^{\frac{2\bar{k}_i}{2\bar{k}-1}},$ 
$(\frac{\bar{b}_{2i}}{r_i})^{\frac{1}{2\bar{k}-1}}$ are integrable within $[0,T].$

\end{remark}
%\section{Discrete-time mean-field-type games}
%This section shows that  direct method can be applied to discrete time mean-field-type games beyond the LQ and variance reduction setting.
\section{Numerical Examples}  \label{sec:mf:num:gv}
In this section, we present some numerical illustrations of  Problem (\ref{gv:problem}) by choosing Gauss-Volterra process with the following kernel
\begin{equation}\begin{array}{ll}
K(t,t')=K_H(t,t')= c_H (t-t')^{H-\frac{1}{2}} \\
+
 c_H (\frac{1}{2}-H)\int_{t'}^t  (z-t')^{H-\frac{3}{2}} (1-(\frac{t'}{z})^{\frac{1}{2}-H}) dz,

\end{array}
\end{equation}
and $H\in (0,1),$ $ c_H=\sqrt{\frac{2H\Gamma( \frac{3}{2}-H)}{\Gamma( \frac{1}{2}+H)\Gamma( 2-2H)}},$ is a normalizing constant, where $\Gamma$ is the gamma function
$$
\Gamma(z)=\int_{0}^{+\infty} e^{-t} t^{z-1} dt,\ \mbox{ Re}(z)>0.
$$
The Gauss-Volterra process with kernel $K_H$ is  a fractional Brownian motion with the  Hurst parameter $H.$
The parameters of the numerical setting are displayed in Table \ref{table:parm:fractional}.

\begin{figure}[htb]
\centering
\begin{tabular}{|c|c|}  
\hline
Numerical setting & Value\\ \hline
Kernel $K$  & $K=K_{0.8} $ \\ \hline
$T$ & $1$\\
Switching $\mathcal{S}$ & $\{s_*, s^*\},\ s_*\neq s^*$\\
$I$ & $2019$\\  \hline
 $(k_i,\bar{k}_i)$ & $ (2,2)$ if $i\leq 2018$ \\ 
  & $ (4,2)$ if $i= 2019$ \\ \hline
  $q_i, \bar{q}_i, r_i,\bar{r}_i$ & $ 1$  if $i\leq  2018$ \\ 
 % $q_{i}, \bar{q}_i, r_i,\bar{r}_i$ & $ 10$  if $i=2008$ \\ \hline
 & $ 100$  if $i=2019$ \\ \hline
$b_1,  b_{2i},\bar{b}_{2i}, \sigma,\sigma_{gv}$ & $ 1$  if $i\leq  2018$  \\ \hline
$ b_{2i},\bar{b}_{2i}$ & $ 10$  if $i= 2019$  \\ \hline
$\mu(t,\theta,s)$ & $ \theta$ \\  
$\Theta$ & $\mathbb{R}^*_{+}$\\
$\nu(d\theta)$ & $ ce^{-5|\theta|}d\theta $\\  \hline
$\sigma_{gv}(t,s^*)$ & $ 1$    \\ \hline
$\sigma_{gv}(t,s_*)$ & $ 10^{-2}$    \\ \hline
 $\tilde{q}_{s^*s_*}$ & $  0.7$\\
$\tilde{q}_{s_*s^*}$ & $ 0.4$\\

\hline % {semi:prop2}
\end{tabular}  \\
\caption{Parameters used in the numerical example.}
\label{table:parm:fractional}
\end{figure}

It is important to notice that  under this setting the problem (\ref{gv:problem}) is not Markov and the cost is not quadratic. From (\ref{semi:prop9}) we know that
the mean-field Nash equilibrium of the mean-field type game (\ref{gv:problem}) under Gauss-Volterra process is 
given by 
\begin{equation}
\begin{array}{ll}
u^*_i=-  \alpha_i^{\frac{1}{3}} (x-\bar{x})-  \bar{\alpha}_i^{\frac{1}{3}} \bar{x},\ \ i\leq 2018\\
u^*_i= -( \frac{\alpha_i}{10})^{\frac{1}{7}} (x-\bar{x})-( \frac{\bar{\alpha}_i}{10})^{\frac{1}{3}} \bar{x},\  i=2019\\
%\mathbb{E}L_i=\mathbb{E}[ \alpha_i(0,s_0) \frac{(x_0-\bar{x}_0)^{2k_i}}{2k_i}  + \bar{\alpha}_i (0,s_0)\frac{\bar{x}_0^{2\bar{k}_i}}{2\bar{k}_i}],\\
%\dot{\alpha}_i  + q_{i}+ 2k_i\alpha_i + \alpha_i k_i(2k_i-1)(1+{\color{blue} \sigma^2_{cogv}} )\\
%+ {\alpha}_i\int_{\theta}[ (1+\theta)^{2k_i}-1 -2k_i\theta   ] \nu(d\theta)\\
%+\sum_{s'} [{\alpha}_i(t,s')- \alpha_i(t,s)]\tilde{q}_{ss'}
%-(2k_i-1)  r_i ( -\frac{b_{2i}\alpha_i}{r_i})^{\frac{2k_i}{2k_i-1}} \\ +2k_i\alpha_i \sum_{j\neq i} b_{2j} ( -\frac{b_{2j}\alpha_j}{r_j})^{\frac{1}{2k_j-1}}=0,\\
%{\alpha}_i(T,s)=1,\\
%%%
% \dot{\bar{\alpha}}_i + \bar{q}_{i} + 2\bar{k}_i\bar{\alpha}_i \bar{b}_1  
%+\sum_{s'} [\bar{\alpha}_i(t,s')- \bar{\alpha}_i(t,s)]\tilde{q}_{ss'} \\
%- (2\bar{k}_i-1) \bar{r}_i ( -\frac{\bar{b}_{2i}\bar{\alpha}_i}{\bar{r}_i})^{\frac{2\bar{k}_i}{2\bar{k}_i-1}} 
%+2\bar{k}_i\bar{\alpha}_i[\sum_{j\neq i}  \bar{b}_{2j} ( -\frac{\bar{b}_{2j}\bar{\alpha}_j}{\bar{r}_j})^{\frac{1}{2\bar{k}_j-1}} ] 
%=0,\\
%\bar{\alpha}_i(T,s)=1,\\
%\sigma_{cogv}^2=\frac{d}{dt}\bigg[\int_0^t   \Bigg\{ K(t'_+,t')\sigma_{gv}(t')+\int_{t'}^{t}\sigma_{gv}(t'')K(t',t'')dt''\Bigg\}^2  dt'\bigg],
 
\end{array}
\end{equation}
Figure \ref{fig:v1} plots (a) a sample path of the  optimal state trajectory starting from $x_0=50,$  (b) the optimal strategies of all decision-makers $i\leq 2018$ and 2019, and (c) sample noises. As expected the state is moving toward zero when the optimal strategies are employed. 
\begin{figure*} \hspace{-2cm}
 \includegraphics[clip, trim=0cm 11cm 0cm 0cm,scale=0.5]{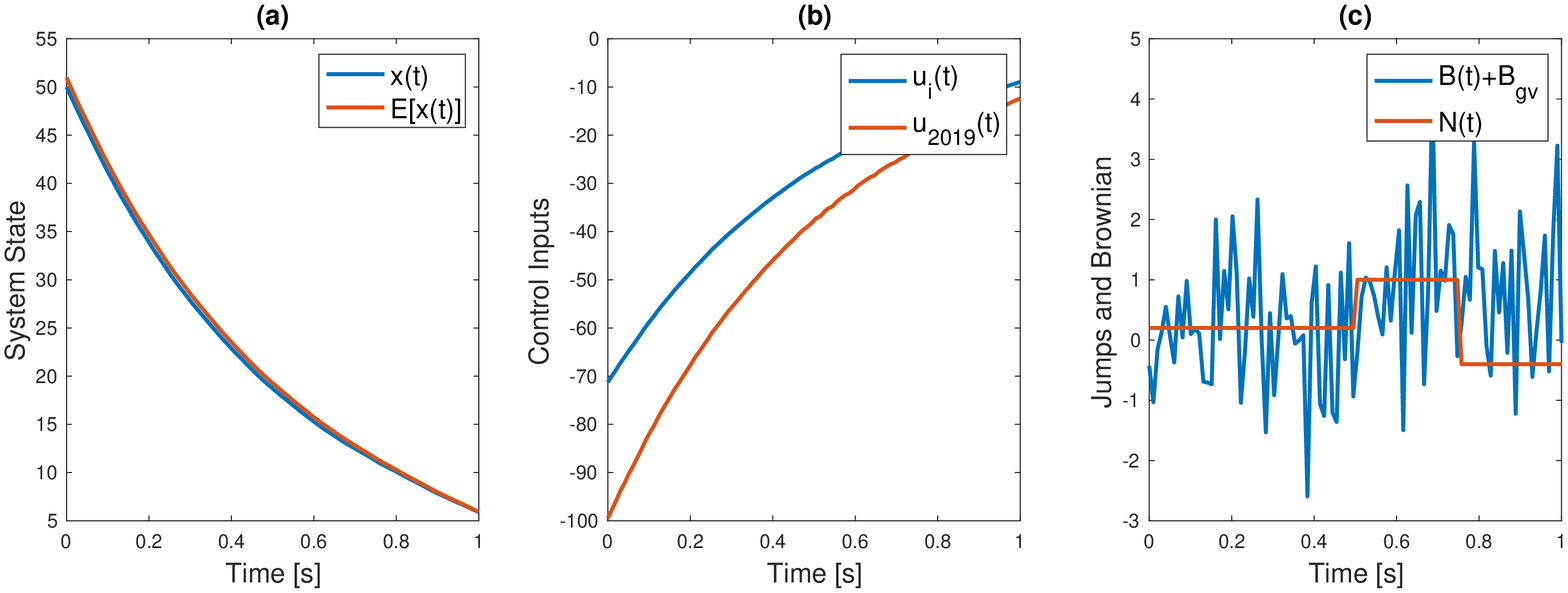}\\ \label{fig:v1}  \caption{Sample path of the optimal state, optimal strategy  under Gauss-Volterra, jump and diffusion. Initial state value is $50.$ }
\end{figure*}
%\section{Discussions and Open issues}
%{\color{blue} to be completed}

\section{Conclusion} %

In this article, we have shown that a mean-field equilibrium can be determined in a semi-explicit way for a broader class of non-linear, non-quadratic game problems with non-linearly distribution-dependent payoffs where the state dynamics is driven by conditional expected values of states, controls, Brownian motions, Gauss-Volterra processes, jump and regime-switching. The method does not require the sophisticated non-elementary extension to backward-forward systems. It does not need PIDEs. It does not need SMPs. It is basic and applies the stochastic integration formula. The use of this simple method may open the accessibility of the tool to a broader audience including beginners and engineers to this emerging field of mean-field-type game theory. Another direct application of the results presented this article is that the explicit solution provides a reference trajectory to the numerical schemes of the corresponding master system beyond the LQ setting. 

\section*{Acknowledgements}
Authors gratefully acknowledge support from U.S. Air Force Office of Scientific Research under grants number FA9550-17-1-0259 , FA9550-12-1-0384 and 
NSF grant DMS 1411412.

\newpage 
%\parpic{\includegraphics[clip, trim=0.4cm 1cm 1cm 2cm,width=0.8in,clip,keepaspectratio]{1_Julian_Barreiro-Gomez.jpg}}

\noindent {\bf Julian Barreiro-Gomez} received his B.S. degree (cum laude) in Electronics Engineering from Universidad Santo Tomas (USTA), Bogot\'{a}, Colombia, in 2011. He received the MSc. degree in Electrical Engineering and the Ph.D. degree in Engineering from Universidad de Los Andes (UAndes), Bogot\'{a}, Colombia, in 2013 and 2017, respectively.  He received the Ph.D. degree (cum laude) in Automatic, Robotics and Computer Vision from the Technical University of Catalonia, Barcelona, Spain, in 2017, and the Best Ph.D. Thesis in Control Engineering 2017 award from the Spanish National Committee of Automatic Control  and Springer. He is currently a Post-Doctoral Associate in the Learning \& Game Theory Laboratory at  New York University Abu Dhabi. His main research interests are evolutionary game dynamics, mean-field-type games, and distributed control and optimization.

%\parpic{\includegraphics[clip, trim=0cm 0cm 0cm 0cm,width=0.8in,clip,keepaspectratio]{duncan_tyrone.jpg}}
\vspace{2cm}
\noindent {\bf Tyrone E. Duncan}  received the B.E.E. degree from Rensselaer Polytechnic Institute, Troy, NY, in 1963 and the M.S. and Ph.D. degrees from Stanford University, Stanford, CA, in 1964 and 1967, respectively. He has held regular positions with the University of Michigan, Ann Arbor (1967-1971), the State University of New York, Stony Brook (1971-1974), and the University of Kansas, Lawrence (1974- present), where he is Professor of Mathematics. He has held visiting positions with the University of California, Berkeley (1969-1970), the University of Bonn, Germany (1978-1979), and Harvard University, Cambridge, MA (1979-1980), and shorter visiting positions at numerous other institutions. Dr. Duncan is a member of the editorial boards of Communications on Stochastic Analysis, and Risk and Decision Analysis and was on the editorial board of SIAM Journal on Control and Optimization (1994-2007) as an Associate Editor and a Corresponding Editor. He is a member of AMS, MAA, and SIAM.

\vspace{2cm}
%\parpic{\includegraphics[clip, trim=0cm 0cm 0cm 0cm,width=0.8in,clip,keepaspectratio]{duncan_bozenna.jpg}}
\noindent {\bf Bozenna Pasik-Duncan}
 received M.S. degree in mathematics from University of Warsaw, and Ph.D. and D.Sc. (Habilitation) degrees from the Warsaw School of Economics, Poland. She is Professor of Mathematics; Courtesy Professor of EECS and AE; Investigator at ITTC; Affiliate Faculty at Center of Computational Biology, and Chancellors Club Teaching Professor at the University of Kansas (KU). She is strong advocate for STEM education and for women in STEM. She is IEEE WIE Chair, founder of CSS Women in Control, founder and faculty advisor of KU Student Chapters of AWM and SIAM, founder and coordinator of KU and CSS Outreach Programs. She has served in many capacities in several societies. Her current service includes Chair of IEEE WIE, member of IEEE CSS and IEEE SSIT Board of Governors, Deputy Chair of IEEE CSS TC on Control Education, Chair of AACC Education Committee, member of IFAC TB, and Award Committees of AWM and MAA. She is General Chair of IFAC ACE 2019 Symposium and member of Organizing Committee of SIAM CT19. She is Associate Editor of several Journals. Her research interests are primarily in stochastic systems and adaptive control and its applications to science and engineering, and in STEM education. She is a recipient of many awards that include IREX Fellowship, NSF Career Advancement Award, Louise Hay Award, Polish Ministry of Higher Education Award, H.O.P.E., Kemper Fellowship, IEEE EAB Meritorious Achievement Award, Service to Kansas and IFAC Outstanding Service Awards. She is Life Fellow of IEEE and Fellow of IFAC, a recipient of the IEEE Third Millennium Medal and IEEE CSS Distinguished Member Award. She is inducted to the KU Women's Hall of Fame.

%\parpic{\includegraphics[clip, trim=7cm 3cm 6.5cm 0cm,width=0.8in,clip,keepaspectratio]{2_Hamidou_Tembine.jpg}}
\vspace{2cm}
\noindent {\bf Hamidou Tembine} received the M.S. degree in applied mathematics from Ecole Polytechnique, Palaiseau, France, in 2006 and the Ph.D. degree in computer science from the University of Avignon, France, in 2009. He holds over 150 scientific publications including magazines, letters, journals, and conferences. He is the author of the book on Distributed Strategic Learning for Engineers (CRC Press, Taylor \& Francis 2012), and coauthor of the book Game Theory and Learning in Wireless Networks (Elsevier Academic Press). 

\end{document}